\documentclass{amsart}
\usepackage{amssymb}

\input epsf

\swapnumbers
\newtheorem{l-theorem}[subsubsection]{Theorem}
\newtheorem{l-proposition}[subsubsection]{Proposition}
\newtheorem{l-lemma}[subsubsection]{Lemma}
\newtheorem{l-corollary}[subsubsection]{Corollary}
\theoremstyle{remark}
\newtheorem{l-remark}[subsubsection]{Remark}
\theoremstyle{definition}
\newtheorem{l-definition}[subsubsection]{Definition}
\numberwithin{equation}{section}

\def\theorem{\begin{l-theorem}}
\def\endtheorem{\end{l-theorem}}
\def\proposition{\begin{l-proposition}}
\def\endproposition{\end{l-proposition}}

\def\corollary{\begin{l-corollary}}
\def\endcorollary{\end{l-corollary}}
\def\remark{\begin{l-remark}}
\def\endremark{\end{l-remark}}
\def\definition{\begin{l-definition}}
\def\enddefinition{\end{l-definition}}

\def\roster{\begin{enumerate}}
\def\endroster{\end{enumerate}}

\let\local\label
\def\loccit#1{$(\ref{#1})$}
\let\ditto\loccit
\def\iref#1{\ref{#1}\loccit}

\def\bigpagebreak{\bigbreak}
\def\medpagebreak{\medbreak}
\def\period{.\spacefactor1000}

{\catcode`\@=11
\gdef\proclaimfont@{\sl}}

\def\proclaimfont{\it}

\def\paragraph{\subsection{}}

\def\ie,{\emph{i.e.},}
\def\eg,{\emph{e.g.},}

\def\cf.{\emph{cf\period}}
\def\via{\emph{via}}
\def\PL{{\sl PL\/}}

\let\AMSbold\mathbf
\let\Bbb\mathbb
\let\Cal\mathcal

\def\bA{\AMSbold A}
\def\bD{\AMSbold D}
\def\bE{\AMSbold E}
\def\bJ{\AMSbold J}

\def\BG#1{\Bbb B_{#1}}
\def\SG#1{\Bbb S_{#1}}
\def\Cp#1{\Bbb P^{#1}}
\def\Rp#1{\Bbb P_{\R}^{#1}}
\def\ls|#1|{\mathopen|#1\mathclose|}
\def\<#1>{\langle#1\rangle}

\def\ind{\operatorname{ind}}
\def\mult{\operatorname{mult}}
\def\omod{\operatorname{mod}}
\def\Sk{\operatorname{Sk}}

\def\CT{\Cal T}
\def\CO{\Cal O}
\def\COP{\CO_{\Cp1}}
\def\Za{Z_{{\rm ar}}}

\let\Ga\alpha
\let\Gb\beta
\let\Gg\gamma
\let\Gd\delta
\let\Gs\sigma

\def\inserthyphen{\ifcat\next a-\fi\ignorespaces}
\let\BLACK\bullet
\let\WHITE\circ
\def\CROSS{\vcenter{\hbox{$\scriptstyle\mathord\times$}}}
\let\STAR*
\def\pblack-{$\BLACK$\futurelet\next\inserthyphen}
\def\pwhite-{$\WHITE$\futurelet\next\inserthyphen}
\def\pcross-{$\CROSS$\futurelet\next\inserthyphen}
\def\pstar-{$\STAR$\futurelet\next\inserthyphen}
\def\black{\protect\pblack}
\def\white{\protect\pwhite}
\def\cross{\protect\pcross}
\def\star{\protect\pstar}
\def\NO#1{\mathord\#_{#1}}
\def\nblack{\NO\BLACK}
\def\nwhite{\NO\WHITE}
\def\ncross{\NO\CROSS}
\def\nstar{\NO\STAR}
\def\insertspace{\ifcat\next a\space\else\ifcat\next\space\space\fi\fi\ignorespaces}
\def\psolid-{solid\futurelet\next\insertspace}
\def\pbold-{bold\futurelet\next\insertspace}
\def\pdotted-{dotted\futurelet\next\insertspace}
\def\solid{\protect\psolid}
\def\bold{\protect\pbold}
\def\dotted{\protect\pdotted}

\def\base{S^2}
\def\Base#1{\base_\Gamma}

\def\bb{\frak b}
\def\bc{\frak c}
\def\bl{\frak l}

\def\bb{m}
\def\bc{c}
\def\bl{l}

\def\ix#1{(#1)}

\def\conj{\operatorname{conj}}
\def\const{\operatorname{const}}
\def\id{\operatorname{id}}
\def\pr{\operatorname{pr}}
\def\R{\Bbb R}
\def\C{\Bbb C}
\def\Z{\Bbb Z}
\def\N{\Bbb N}
\def\Q{\Bbb Q}
\let\sminus\smallsetminus

\let\ge\geqslant
\let\le\leqslant

\def\rom#1{\/{\rm#1}}

\begin{document}

\author{Alex Degtyarev}

\title{Zariski $k$-plets via dessins d'enfants}

\address{Department of Mathematics,
Bilkent University,
06800 Ankara, Turkey}

\email{degt@fen.bilkent.edu.tr}

\begin{abstract}
We construct exponentially large collections of pairwise
distinct equisingular deformation families of irreducible plane
curves sharing the same sets of singularities. The fundamental
groups of all curves constructed are abelian.
\end{abstract}

\keywords{Zariski pair, trigonal curve, \emph{dessin d'enfants},
braid monodromy}

\subjclass[2000]{
Primary: 14H50; 
Secondary: 14H30, 
14D05 
}

\maketitle

\section{Introduction}\label{S.intro}

\subsection{Motivation and principal results}\label{results}
Throughout this paper, the \emph{type} of a singular point
is its embedded \emph{piecewise linear} type, and
equisingular deformations of curves
in surfaces are understood in the \emph{piecewise linear} sense,
\ie, the \PL-type of each singular point should be preserved
during the deformation. This convention is essential as some
of the curves considered have
non-simple singularities.

Recall that a \emph{Zariski $k$-plet} is a collection
$C_1,\ldots,C_k$ of plane curves, all of the same degree~$m$, such
that
\roster
\item\local{comb.data}
all curves have the same \emph{combinatorial data}
(see~\cite{Artal.comb}
for the definition; for irreducible curves, this means the
set of types of singular points), and
\item\local{def.type}
the curves are pairwise not equisingular deformation equivalent.
\endroster
Note that Condition~\loccit{def.type} in the definition differs
from paper to paper, the most common being the requirement that
the pairs $(\Cp2,C_i)$ (or complements
$\Cp2\sminus C_i$) should not be homeomorphic.
In this paper, we choose equisingular deformation equivalence,
\ie, being in the same component of the moduli space, as it is the
strongest topologically meaningful `global' equivalence relation.
In any case, the construction of topologically distinguishable
Zariski $k$-plets usually starts with finding curves
satisfying~\loccit{def.type} above.

Historically, the first example of Zariski pairs was found by
O.~Zariski~\cite{Zariski.group}, \cite{Zariski}, who constructed a
pair of irreducible sextics~$C_1$, $C_2$, with six cusps each,
which differ by the fundamental groups $\pi_1(\Cp2\sminus C_i)$.
Since then, a great number of other examples has been found.
Citing recent results only, one
can mention a large series of papers by
E.~Artal Bartolo, J.~Carmona Ruber, J.~I.~Cogolludo Agust\'{\i}n,
and H.~Tokunaga (see~\cite{Artal.comb}, \cite{Artal.k-plets} and
more recent papers~\cite{Artal.Trends}--\cite{Artal.KT} for
further references), A.~Degtyarev~\cite{poly}, \cite{JAG},
\cite{degt.Oka}
(paper~\cite{poly} deals with a direct generalization of
Zariski's example: pairs of sextics distinguished by their
Alexander polynomial), C.~Eyral and M.~Oka~\cite{EO1}, \cite{EO2},
\cite{Oka},
G.-M.~Greuel, C.~Lossen, and E.~Shustin~\cite{Shustin}
(Zariski pairs with abelian fundamental groups),
Vik.~S.~Kulikov~\cite{Kulikov.pairs},
A.~\"Ozg\"uner~\cite{Aysegul} (a complete list of Zariski pairs
of irreducible sextics that are distinguished by their
Alexander polynomial),
I.~Shimada~\cite{Shimada1}--\cite{Shimada3} (a complete list of
Zariski pairs of sextics with the maximal total Milnor number
$\mu=19$, as well as a list of arithmetic Zariski pairs of
sextics),
and A.~M.~Uluda\u{g}~\cite{Uludag}.
The amount of literature on the subject definitely calls for a
comprehensive survey!

With very few exceptions, the examples found in the literature are
those of Zariski \emph{pairs} or \emph{triples}. To my knowledge,
the largest known Zariski $k$-plets are those constructed
in Artal Bartolo, Tokunaga~\cite{Artal.k-plets}:
for each integer $m\ge6$, there is a
collection of $([m/2]-1)$ reducible curves of degree~$m$ sharing
the same combinatorial data. The principal result of this paper is
the following Theorem~\ref{th.main}, which states that the size
of Zariski $k$-plets can grow exponentially with the degree.
(Theorem~\ref{th.reducible} below gives a slightly better count
for reducible curves.)

\theorem\label{th.main}
For each integer $m\ge8$, there is a set of singularities shared by
$$
Z(m)=\frac1k\binom{2k-2}{k-1}\binom{k}{[k/2]}\binom{[k/2]}\epsilon
$$
pairwise distinct equisingular deformation families of
irreducible plane curves~$C_i$ of
degree~$m$, where
$k=[(m-2)/2]$ and $\epsilon=m-2k-2\in\{0,1\}$.
The fundamental groups of all curves~$C_i$ are abelian\rom: one has
$\pi_1(\Cp2\smallsetminus C_i)=\Z_m$.
\endtheorem

Recall that a \emph{real structure} on a complex surface~$X$ is an
anti-holomorphic involution $\conj\colon X\to X$. A curve $C\subset X$
is called \emph{real} (with respect to~$\conj$) if $\conj(C)=C$,
and a deformation $C_t$, $\ls|t|\le1$, is called \emph{real} if
$C_{\bar t}=\conj C_t$. Up to projective equivalence, there is a
unique real structure on~$\Cp2$; in appropriate coordinates it is
given by $(z_0:z_1:z_2)\mapsto(\bar z_0:\bar z_1:\bar z_2)$.

For completeness, we enumerate the families containing real
curves.

\theorem\label{th.real}
If $m=8t+2$ for some $t\in\Z$, then $Z(4t+2)$ of the
families
given
by Theorem~\ref{th.main}
contain real curves \rom(with respect to
some real structure in~$\Cp2$\rom). All other curves \rom(and all
curves for other values of~$m$\rom) split into pairs of disjoint
complex conjugate equisingular deformation families.
\endtheorem

\theorem\label{th.reducible}
For each integer $m\ge8$, there is a set of combinatorial data
shared by
$$
R(m)=\frac1{m-5}\binom{2m-12}{m-6}
$$
pairwise distinct equisingular deformation families of
plane curves~$C_i$ of
degree $m$ \rom(each curve splitting into an irreducible
component of degree $(m-1)$ and a line\rom).
The fundamental groups of all curves~$C_i$ are abelian\rom: one has
$\pi_1(\Cp2\smallsetminus C_i)=\Z$.

If $m=2t+1$ is odd, then $R(t+3)$ of the families above contain
real curves \rom(with respect to
some real structure in~$\Cp2$\rom). All other curves \rom(and all
curves for $m$~even\rom) split into pairs of disjoint
complex conjugate equisingular deformation families.
\endtheorem

Theorems~\ref{th.main}--\ref{th.reducible} are proved in
Sections~\ref{proof.main}--\ref{proof.reducible}, respectively.

It is easy to see that the counts~$Z(m)$ and~$R(m)$ given by the
theorems
grow faster
than $a^{3m/2}$ and~$a^{2m}$, respectively, for any $a<2$. A few values
of~$Z$ and~$R$ are listed in the table below.

\toks0{m}
\toks2{Z(m)}
\toks4{R(m)}
\def\add[#1,#2,#3]{\toks0\expandafter{\the\toks0&\hss$#1$\hss}
 \toks2\expandafter{\the\toks2&\hss$#2$\hss}
 \toks4\expandafter{\the\toks4&\hss$#3$\hss}}
\add[8, 6, 2]
\add[9, 6, 5]
\add[10, 30, 14]
\add[11, 60, 42]
\add[12, 140, 132]
\add[13, 280, 429]
\add[14, 840, 1430]
\add[\!\!\ldots\!\!, \!\!\ldots\!\!, \!\!\ldots\!\!]
\add[20, 2\cdot10^5, 3\cdot10^6]
\add[40, 4\cdot10^{13}, 8\cdot10^{17}]
\add[80, 1\cdot10^{31}, 3\cdot10^{41}]

\bigpagebreak\centerline{\vbox{%
\halign{\vrule depth3pt height10pt\quad\hss$#$\hss\quad\vrule&&
 \quad#\cr
\noalign{\hrule}
\the\toks0\quad\vrule\cr
\noalign{\hrule}
\the\toks2\quad\vrule\cr
\noalign{\hrule}
\the\toks4\quad\vrule\cr
\noalign{\hrule}}}}
\medpagebreak

Note that we are not trying to set a record here; probably, there
are much larger collections of curves constituting Zariski
$k$-plets. The principal emphasis of this paper is the fact that
Zariski $k$-plets can be exponentially large.

\subsection{Other results and tools}\label{tools}
The curves given by Theorems~\ref{th.main} and~\ref{th.reducible}
are plane curves of degree~$m$ with a singular point of
multiplicity $(m-3)$. (In a sense, this is the first nontrivial
case, as curves with a singular point of multiplicity $(m-2)$ or
$(m-1)$ do not produce Zariski pairs, see~\cite{quintics}.) When
the singular point is blown up, the proper transform of the curve
becomes a (generalized) trigonal curve in a rational ruled
surface. We explain this relation in Section~\ref{S.models}, and
the bulk of the paper deals with trigonal curves, whose theory is
rather parallel to Kodaira's theory of Jacobian
elliptic fibrations.

A trigonal curve can be characterized by its functional
$j$-invariant, which is a rational function $j\colon \Cp1\to\Cp1$, so
that the singular fibers of the curve are encoded in terms of
the pull-back $j^{-1}\{0,1,\infty\}$ (see Table~\ref{table.j}). To
study the $j$-invariants, we follow S.~Orevkov's
approach~\cite{Orevkov}, \cite{Orevkov-private} (see
also~\cite{DIK.elliptic}) and use a modified version of
Grothendieck's \emph{dessins d'enfants}, see
Section~\ref{S.dessins},  reducing the classification of trigonal
curves with prescribed combinatorial type of singular fibers to a
graph theoretical problem. The resulting problem is rather
difficult, as the graphs are allowed to undergo a number of
modifications (see~\ref{dessin.equivalence})
caused by the fact that $j$ may have critical values
other than~$0$, $1$, or~$\infty$. To avoid this difficulty, we
concentrate on a special case of the so called \emph{maximal}
curves, see~\ref{def.maximal}, which can be characterized as
trigonal curves not admitting any further degeneration
(Proposition~\ref{no.degeneration}); the classification of maximal
curves reduces to the enumeration of connected
planar maps with vertices of
valency~$\le3$, see Theorem~\ref{B->Sk}. We exploit this relation
and use oriented rooted binary trees to produce large Zariski
$k$-plets of
trigonal curves, see Proposition~\ref{C.curves} and a slight
modification in Proposition~\ref{alt.curves}.

It is worth mentioning that the curves given by
Propositions~\ref{C.curves} and~\ref{alt.curves} are defined over
algebraic number fields (like all maximal curves), and in
Theorem~\ref{alt.main} we use this fact to construct a slightly
smaller, but still exponentially large, Zariski $k$-plet of plain
curves with discrete moduli space. All these examples seem to be
good candidates for exponentially large \emph{arithmetic Zariski
$k$-plets} (in rational ruled surfaces and in the plane)
in the sense of Shimada~\cite{Shimada1},
\cite{Shimada2}.

An important question that remains open is whether the curves
constituting various Zariski $k$-plets constructed in the paper
can be distinguished topologically. As a first step in this
direction, we calculate the braid monodromy of the trigonal
curves, see~\ref{C.monodromy}. (For the relation between the braid
monodromy and the topology of the curve, see
Orevkov~\cite{Orevkov-private},
Vik.~S.~Kulikov and M.~Teicher~\cite{Kulikov},
or Carmona Ruber~\cite{Carmona}.) In~\ref{monodromy}, we give a
general description of the braid monodromy of a trigonal curve in
terms of its dessin; it covers all maximal curves with the
exception of four explicitly described series. As a simple
application, we obtain a criterion of reducibility of a maximal
trigonal curve in terms of its skeleton, see
Corollary~\ref{reducible.curve}.

As another direct application of the construction, we produce
exponentially large Zariski $k$-plets of Jacobian elliptic
surfaces, see~\ref{surfaces}. (Here, by a Zariski $k$-plet we mean
a collection of not fiberwise deformation equivalent surfaces
sharing the same combinatorial type of singular fibers.) The
series given by Theorem~\ref{alt.surfaces} are related to positive
definite lattices of large rank; this gives one hope to
distinguish the surfaces, and hence their branch loci,
topologically.

\subsection{Contents of the paper}\label{contents}
In Section~\ref{S.models}, we introduce trigonal curves in
rational ruled surfaces and discuss their relation to plane
curves with a singular point of multiplicity $\text{degree}-3$.
Section~\ref{S.j} reminds the basic properties of the
$j$-invariant of a trigonal curve, and Section~\ref{S.dessins}
introduces the dessin of a trigonal
curve and the skeleton of a maximal
curve. In Section~\ref{S.group}, we prove a few technical
statements on the fundamental group of a generalized trigonal
curve. Section~\ref{S.monodromy} deals with the braid monodromy.
The principal results of the paper,
Theorems~\ref{th.main}--\ref{th.reducible}, are proved in
Section~\ref{S.construction}.
Finally, in
Section~\ref{S.applications}, we discuss a few modifications of
the construction and state a few open problems.

\section{Trigonal models}\label{S.models}

\subsection{Hirzebruch surfaces}\label{Hirzebruch}
Recall that the \emph{Hirzebruch surface} $\Sigma_k$, $k\ge0$, is
a rational geometrically ruled surface with a section~$E$ of
self-intersection~$-k$. If $k>0$, the ruling is unique and there
is a unique section~$E$ of self-intersection~$-k$; it is called
the \emph{exceptional section}.
In the
exceptional case $k=0$, the surface $\Sigma_0=\Cp1\times\Cp1$
admits two rulings, and we choose and fix one of them;
any fiber of the other ruling can be chosen for the
exceptional section. The fibers
of the ruling are referred to as the fibers of~$\Sigma_k$.
The semigroup of classes of effective divisors
on
~$\Sigma_k$
is generated by the classes of the exceptional
section~$E$ and a fiber~$F$; one has
$E^2=-k$, $F^2=0$, and $E\cdot F=1$.

An \emph{elementary transformation} of a Hirzebruch
surface~$\Sigma_k$ is the birational transformation consisting in
blowing up a point $O\in\Sigma_k$ and blowing down the proper
transform of the fiber through~$O$. If the blow-up center~$O$ does
(respectively, does not) belong to the exceptional section
$E\subset\Sigma_k$, the result of the elementary transformation is
the Hirzebruch surface~$\Sigma_{k+1}$ (respectively,
$\Sigma_{k-1}$).

\subsection{Trigonal curves}\label{trigonal.curves}
A \emph{generalized trigonal curve} on a
Hirzebruch surface~$\Sigma_k$ is a reduced curve
not containing the exceptional section~$E$ and
intersecting each generic fiber at three
points. Note that a generalized trigonal curve $B\subset\Sigma_k$
may contain fibers of~$\Sigma_k$ as components; we will call them
the \emph{linear components} of~$B$.

A \emph{singular fiber} of a generalized trigonal curve
$B\subset\Sigma_k$ is a fiber~$F$ of~$\Sigma_k$ that is not
transversal to the union $B\cup E$. Thus, $F$ is either a linear
component of~$B$, or the fiber through a point of intersection
of~$B$ and~$E$, or the fiber over a critical value of the
restriction to~$B$ of the projection $\Sigma_k\to\Cp1$.

A \emph{trigonal curve} is a generalized trigonal curve
disjoint from
the exceptional section. (In particular, trigonal curves have no
linear components.)
For a trigonal curve $B\subset\Sigma_k$,
one has $\ls|B|=\ls|3E+3kF|$; conversely, any curve
$B\in\ls|3E+3kF|$ not containing~$E$ as a component is a trigonal
curve.

Let~$F$ be a singular fiber of a trigonal curve~$B$.
If $B$ has at most simple singular points
on~$F$ and $F$ is not a component of~$B$,
then locally $B\cup E$ is the branch locus of a Jacobian
elliptic surface~$X$, and the pull-back of~$F$ is a singular fiber
of~$X$. In this case, we use the standard
notation for singular elliptic fibers
(referring to the extended Dynkin diagrams) to describe the type
of~$F$. Otherwise, $B$ has a singular
point of type~$\bJ_{k,p}$ or~$\bE_{6k+\epsilon}$, see~\cite{AVG}
for the notation,
and we use the notation $\tilde\bJ_{k,p}$
and~$\tilde\bE_{6k+\epsilon}$, respectively, to describe the type
of~$F$.

\remark
We will not attempt to give a formal definition of the type of a
singular fiber~$F$ of a trigonal curve~$B$. One can understand it
as the topological type of the boundary singularity $(B,F)$,
see~\cite{AVG} for details. As a result of the classification, one
can conclude that this type is determined by whether $F$ is a
component of~$B$ and (the conjugacy class of)
the braid monodromy about~$F$, see
Section~\ref{S.monodromy} below for the definition. Alternatively,
if $F$ is not a component of~$B$ and $B$ has at worst simple
singularities on~$F$ (which is always the case in this paper),
then the type of~$F$ is determined by Kodaira's type of the singular
fiber of the Jacobian elliptic surface
ramified at $B\cup F$, see above.
\endremark

Any generalized trigonal curve~$B$ without linear
components can be converted to a trigonal curve by a sequence of
elementary transformations, at each step blowing up a point of
intersection of~$B$ and the exceptional section and blowing down
the corresponding fiber.

\subsection{Simplified models}\label{simplified.models}
Let $\Sigma'$ be a Hirzebruch surface, and let $\Sigma''$ be
obtained from~$\Sigma'$ by an elementary transformation.
Denote by $O'\in\Sigma'$ and $O''\in\Sigma''$ the blow-up centers
of the transformation and its inverse, respectively, and let
$F'\subset\Sigma'$ and $F''\subset\Sigma''$ be the fibers
through~$O'$ and~$O''$, respectively.
The \emph{transform} $B''\subset\Sigma''$ of a generalized
trigonal curve $B'\subset\Sigma'$ is defined as follows: if $B'$
does not (respectively, does) contain~$F'$ as a linear component,
then $B''$ is the proper transform of~$B'$ (respectively, the
union of the proper transform and fiber~$F''$). In the above
notation, there is an obvious diffeomorphism
\begin{equation}
\Sigma'\sminus(B'\cup E'\cup F')\cong
 \Sigma''\sminus(B''\cup E''\cup F''),\label{elementary.diffeo}
\end{equation}
where $E'\subset\Sigma'$ and $E''\subset\Sigma''$ are the
exceptional sections.

A trigonal curve $B\subset\Sigma_k$ is called \emph{simplified} if
all its singular points are double, \ie, those of type~$\bA_p$.
Clearly, each trigonal curve has a unique \emph{simplified model}
$\bar B\subset\Sigma_l$, which is obtained from~$B$ by a series of
elementary transformations: one blows up a triple point of the curve
and blows down the corresponding fiber, repeating this process
until there are no triple points left.

\subsection{Deformations}\label{deformations}
Let $B\subset\Sigma_k$ be a generalized trigonal curve and
$E\subset\Sigma_k$ the exceptional section.
We define a \emph{fiberwise deformation} of~$B$ as an equisingular
deformation (path in the space of curves) preserving the
topological types of all singular fibers. Alternatively, a
fiberwise deformation can be defined as an equisingular
deformation of the curve
$B\cup E\cup\text{(all singular fibers of~$B$)}$.

A \emph{degeneration} of a generalized trigonal curve~$B$ is a
family $B_t$, $\ls|t|\le1$, of generalized trigonal curves such
that $B=B_1$ and the restriction of~$B_t$ to the annulus
$0<\ls|t|\le1$ is a fiberwise deformation. A degeneration is
called \emph{nontrivial} if $B_0$ is not fiberwise deformation
equivalent to~$B$.

Let $B_k\subset\Sigma_k$ and $B_{k+1}\subset\Sigma_{k+1}$ be two
generalized trigonal curves related by an elementary
transformation, and let $E_i\subset\Sigma_i$, $i=k,k+1$, be the
respective exceptional sections. In general, it is \emph{not} true
that an equisingular deformation of~$B_k$ or $B_k\cup E_k$ is
necessarily followed by an equisingular deformation of~$B_{k+1}$
(respectively, $B_{k+1}\cup E_{k+1}$) or \emph{vice versa}: it may
happen that a singular fiber splits into two and this operation
affects the topology of one of the curves without affecting the
topology of the other.
However, it obviously \emph{is} true that the
fiberwise deformations
of~$B_k$ are in a natural one-to-one correspondence with the
fiberwise deformations of~$B_{k+1}$. A precise statement relating
deformations of~$B_k$ and~$B_{k+1}$ would require simple but
tedious analysis of a number of types of singular fibers.
Instead of attempting to study this problem in full generality
(which becomes even more involved if the two curves are related by a
series of elementary transformations), we just make sure that,
in the examples considered in this paper
(see~\ref{proof.main}, \ref{proof.reducible},
and~\ref{alt.main}),
generic equisingular deformations of each curve~$B\cup E$
are fiberwise. (In~\ref{proof.reducible}, a linear component is
added to the curve for this purpose.) In more details this issue
is addressed in~\ref{rem.deformations}.

\subsection{The trigonal model of a plane curve}\label{trigonal.model}
Let $C\subset\Cp2$ be a reduced curve, $\deg C=m$,
and let~$O$ be a distinguished singular point of~$C$ of
multiplicity~$(m-3)$. (Such a point is unique whenever
$m>7$ or $m>6$ and $C$ is irreducible.) By a \emph{linear
component} of~$C$ we mean a component of degree~$1$ passing
through~$O$.

Blow~$O$ up and denote the result by~$Y_1$; it is a Hirzebruch
surface~$\Sigma_1$, and the proper transform $\tilde C=B_1\subset Y_1$
of~$C$ is a generalized trigonal curve. Clearly, the combinatorial
type of~$C$ determines and is determined by that of $B_1\cup E_1$,
the type of~$O$ itself being recovered from the singularities
of~$B_1\cup E_1$ located in the exceptional section~$E_1$.
Furthermore, equisingular deformations of the pair $(C,O)$ are in
a one-to-one correspondence with equisingular deformations of
$B_1\cup E_1$.

Let~$B_1'$ be the curve obtained from~$B_1$ by removing its linear
components. As in~\ref{trigonal.curves}, one can apply a sequence
of elementary transformations to get a sequence of curves
$B_i,B_i'\subset Y_i\cong\Sigma_i$, $i=1,\ldots,k$, so that $B_k'$
is a true trigonal curve. (Here, $B_{i+1}$ is the transform
of~$B_i$, and $B_{i+1}'$ is obtained from~$B_{i+1}$ by removing
its linear components. In other words, we pass to the trigonal
model of~$B_1'$ while keeping track of the linear components
of~$C$.) The curve~$B_k'$ is called the \emph{trigonal model}
of~$C$. Finally, passing from~$B_k'$ to its simplified model
$B'\subset Y\cong\Sigma_l$, one obtains the \emph{simplified
trigonal model}~$B'$ of~$C$.

\section{The $j$-invariant}\label{S.j}

The contents of this section is a translation to the language of
trigonal curves of certain well known notions and facts about
elliptic surfaces; for more details
we refer to the excellent founding paper
by K.~Kodaira~\cite{Kodaira} or to more recent
monographs~\cite{FM} and~\cite{BPV}. In the theory of elliptic
surfaces, trigonal curves (in the sense of this paper)
arise as the
branch loci of the Weierstra\ss{} models of
Jacobian
elliptic surfaces
over a rational base.
These curves have at most simple singularities
and belong to even Hirzebruch surfaces $\Sigma_{2s}$.
However, most notions and statements extend,
more or less directly, to
trigonal curves in odd Hirzebruch surfaces
$\Sigma_{2s+1}$.

\subsection{Weierstra{\ss} equation}\label{Weierstrass}
Let $\Sigma_k\to\Cp1$ be a Hirzebruch surface.
Any trigonal curve $B\subset\Sigma_k$
can be given by a \emph{Weierstra{\ss} equation};
in appropriate affine charts it has the form
$$
x^3+g_2x+g_3=0,
$$
where $g_2$ and~$g_3$ are certain sections of~$\COP(2k)$
and~$\COP(3k)$,
respectively, and $x$ is a coordinate such that $x=0$ is the zero
section and $x=\infty$ is the exceptional section $E\subset\Sigma_k$.
The
sections~$g_2$, $g_3$ are determined by the curve uniquely up to
the transformation
\begin{equation}
(g_2,g_3)\mapsto (t^2g_2, t^3g_3),\quad
t\in\C^*.
\label{trigequiv}
\end{equation}

The following statement is straightforward.

\proposition\label{simplified}
A trigonal curve~$B$ as
in~\ref{Weierstrass}
is simplified if and only if there is no point $z\in\Cp1$ which is
a root of~$g_2$ of multiplicity~$\ge2$ and a root of~$g_3$ of
multiplicity~$\ge3$.
\qed
\endproposition

\paragraph\label{j-invariant}
The \emph{\rom(functional\rom) $j$-invariant} of a trigonal curve
$B\subset\Sigma_k$
is the meromorphic function $j=j_B\colon \Cp1\to\Cp1$ given by
$$
j=\frac{4g_2^3}\Delta,\quad \Delta=4g_2^3+27g_3^2,
$$
where $g_2$ and~$g_3$ are the coefficients of the Weierstra{\ss}
equation of~$B$, see~\ref{Weierstrass}.
Here, the domain of~$j$ is the base of the ruling $\Sigma_k\to\Cp1$,
whereas its range is the standard projective line
$\Cp1=\C^1\cup\{\infty\}$. If the fiber~$F_z$ over
$z\in\Cp1$ is nonsingular, then the value $j(z)$ is the usual
$j$-invariant (divided by the magic number $1728=12^3$) of the
quadruple of points cut on~$F_z$ by the union $B\cup E$ (or, in
more conventional terms, the $j$-invariant of the elliptic curve
that is the double of $F_z\cong\Cp1$ ramified at the four points
above). The values of~$j$ at the finitely many remaining points
corresponding to the singular fibers of~$B$ are obtained by analytic
continuation.

Since $j_B$ is defined \via\ affine charts and analytic
continuation, it is obviously invariant under elementary
transformations. In particular, the notion of $j$-invariant can be
extended to generalized trigonal curves (by ignoring the linear
components and passing to a trigonal model), and
the $j$-invariant of a trigonal
curve~$B$ is the same as that of the simplified model of~$B$.

\paragraph
The $j$-invariant $j_B\colon \Cp1\to\Cp1$ has three `special' values: $0$,
$1$, and~$\infty$. The correspondence between the type of a
fiber~$F_z$, see remark in Section~\ref{trigonal.curves},
and the value~$j(z)$ (and the ramification index $\ind_zj$
of~$j$ at~$z$) is shown in Table~\ref{table.j}. (We confine
ourselves to the curves with at worst simple singular points. In
fact,
in view of the invariance of the $j$-invariant
under elementary
transformations, it would suffice to consider type~$\tilde\bA$
singular fibers only. For the reader's convenience, we also cite
Kodaira's notation for the types of singular elliptic fibers, \cf.
Section~\ref{trigonal.curves}.)
If $B$ is a curve in~$\Sigma_k$, the maximal degree
of~$j_B$ is $6k$. However, $\deg j_B$ drops if $B$ has triple
singular points or type~$\tilde\bA_0^{**}$, $\tilde\bA_1^*$,
or~$\tilde\bA_2^*$ singular fibers, see $\Delta\deg j$ in
Table~\ref{table.j}. It is worth mentioning that the $j$-invariant
of a generic trigonal curve is highly \emph{non}-generic, as it
takes values~$0$ and~$1$ with multiplicities~$3$ and~$2$
respectively (see Comments to Table~\ref{table.j}); conversely, a
generic function $j\colon \Cp1\to\Cp1$ would arise as the $j$-invariant
of a trigonal curve with a large number of type~$\tilde\bA_0^{**}$
and~$\tilde\bA_1^*$ singular fibers.

\begin{table}[ht]
\caption{The values~$j(z)$ at singular fibers~$F_z$}\label{table.j}
\def\tabrule{\vrule height11pt depth4pt}\def\neg{\!\!\!}
\def\I{{\rm I}}\def\II{{\rm II}}\def\III{{\rm III}}\def\IV{{\rm IV}}
\centerline{\vbox{\offinterlineskip
\halign{\tabrule\quad\hss#\hss\quad\vrule&&\quad\hss$#$\hss\quad\vrule\cr
\noalign{\hrule}
\multispan2\tabrule\hss Type of~$F_z$\hss\vrule&
 \neg j(z)\neg&\ind_zj&\Delta\deg j&\mult F_z\cr
\noalign{\hrule}
$\neg\tilde\bA_p$ ($\tilde\bD_{p+5}$), $p\ge1$\neg&
 \neg\I_{p+1}\ (\I^*_{p+1})\neg&
 \infty&p+1&0\ (-6)&\neg p+1\ (p+7)\neg\cr
\noalign{\hrule}
$\tilde\bA_0^*$ ($\tilde\bD_5$)&\I_1\ (\I^*_1)&
 \infty&1&0\ (-6)&1\ (7)\cr
\noalign{\hrule}
$\tilde\bA_0^{**}$ ($\tilde\bE_6$)&\II\ (\II^*)&
 0&\neg1\bmod3\neg&-2\ (-8)&2\ (8)\cr
\noalign{\hrule}
$\tilde\bA_1^*$ ($\tilde\bE_7$)&\III\ (\III^*)&
 1&\neg1\bmod2\neg&-3\ (-9)&3\ (9)\cr
\noalign{\hrule}
$\tilde\bA_2^*$ ($\tilde\bE_8$)&\IV\ (\IV^*)&
 0&\neg2\bmod3\neg&\neg-4\ (-10)\neg&4\ (10)\cr
\noalign{\hrule}\crcr}}}
\bigskip\vbox{\parindent\normalparindent
\rightskip\normalparindent \leftskip\normalparindent \noindent\small
{\bf Comments.\enspace}
Fibers of type~$\tilde\bA_0$ (Kodaira's
${\rm I}_0$) are not singular.
For a
nonsingular fiber~$F_z$ with complex multiplication of order~$2$
(respectively,~$3$) one has $j(z)=1$ and
$\ind_zj=0\bmod2$ (respectively, $j(z)=0$ and
$\ind_zj=0\bmod3$).
Singular fibers of type~$\tilde\bD_4$
(Kodaira's ${\rm I}_0^*$)
are
not detected by the $j$-invariant, except that each such fiber
decreases the degree of~$j$ by~$6$.
The multiplicity $\mult F_z$
is the number of
simplest (\ie, type~$\tilde\bA_0^*$) singular fibers resulting from a
generic perturbation of~$F_z$.
\vskip-\bigskipamount}
\end{table}

\subsection{Isotrivial curves}\label{isotrivial}
A trigonal curve $B\subset\Sigma_k$ is called \emph{isotrivial} if
$j_B=\const$. All simplified isotrivial curves can easily be
classified.
\roster
\item\local{j=0}
If $j_B\equiv0$, then $g_2\equiv0$ and $g_3$ is a section of
$\COP(3k)$ whose all roots are simple or double,
see Proposition~\ref{simplified}. The singular fibers of~$B$ are of
type~$\tilde\bA_0^{**}$ (over the simple roots of~$g_3$)
or~$\tilde\bA_2^*$ (over the double roots of~$g_3$).
\item\local{j=1}
If $j_B\equiv1$, then $g_3\equiv0$ and $g_2$ is a section of
$\COP(2k)$ with simple roots only, see
Proposition~\ref{simplified}.
All singular fibers of~$B$ are of
type~$\tilde\bA_1^*$ (over the roots of~$g_2$).
\item\local{j=const}
If $j_B=\const\ne0,1$, then $g_2^3\equiv\lambda g_2^2$ for some
$\lambda\in\C^*$; in view of
Proposition~\ref{simplified}, this implies that
$k=0$ and $g_2,g_3=\const$, \ie, $B$ is a union of disjoint
sections of~$\Sigma_0$. (In particular, $B$ has no singular
fibers.)
\endroster
Note that an isotrivial trigonal curve cannot be fiberwise
deformation equivalent to a non-isotrivial one, as a non-constant
$j$-invariant~$j_B$ would take value~$\infty$ and
hence the curve would have a singular fiber of type~$\tilde\bA_0^*$
or~$\tilde\bA_p$, $p>0$, see Table~\ref{table.j}.

\proposition\label{j->B}
Any non-constant meromorphic function $j\colon \Cp1\to\Cp1$ is the
$j$-invariant of a certain simplified trigonal curve
$B\subset\Sigma_k$\rom;
the latter is unique up to the change of coordinates given
by~\eqref{trigequiv}.
\endproposition

\proof
For simplicity, restrict all
functions/sections to an affine portion $\C^1\subset\Cp1$, which
we assume to contain all pull-backs $j^{-1}(0)$ and $j^{-1}(1)$.
Represent the function $l=j/(1-j)$ by an irreducible fraction
$p/q$. Since $l(\infty)\ne0,1$, one has $\deg p=\deg q$.
For each root~$a$ of~$p$ of multiplicity $1\bmod3$
(respectively, $2\bmod3$), multiply both~$p$ and~$q$ by $(z-a)^2$
(respectively, $(z-a)^4$), and for each root~$b$ of~$q$ of
multiplicity $1\bmod2$, multiply both~$p$ and~$q$ by $(z-b)^3$.
In the resulting representation $l=\bar p/\bar q$,
the multiplicity of each root of~$\bar p$ (respectively,~$\bar q$) is
divisible by~$3$ (respectively,~$2$), and $\bar p$ and~$\bar q$ have no
common roots of multiplicity~$\ge6$. Hence, one has
$\bar p=4g_2^3$ and $\bar q=27g_3^2$
for some polynomials $g_2$, $g_3$ satisfying the condition
in Proposition~\ref{simplified}, and the function
$j=l/(l+1)=\bar p/(\bar p+\bar q)$
is the $j$-invariant of the
simplified trigonal curve $B\subset\Sigma_k$
given by the Weierstra\ss{} equation with coefficients~$g_2$, $g_3$,
where $k=\frac16\deg\bar p=\frac16\deg\bar q$.
Clearly, the polynomials~$g_2$,
~$g_3$ as above are defined by~$l$ uniquely up to the
transformation given by~\eqref{trigequiv}.
\endproof

\proposition\label{j.deformations}
A fiberwise deformation of a non-isotrivial
trigonal curve~$B$
results in a deformation of its $j$-invariant $j=j_B\colon \Cp1\to\Cp1$
with the following properties\rom:
\roster
\item\local{degree}
the degree of the map $j\colon \Cp1\to\Cp1$ remains constant\rom;
\item\local{poles}
distinct poles of~$j$ remain distinct, and their multiplicities
remain constant\rom;
\item\local{mod3}
the multiplicity of each root of~$j$ remains
constant~$\omod3$\rom;
\item\local{mod2}
the multiplicity of each root of~$j-1$ remains
constant~$\omod2$.
\endroster
Conversely, any deformation of nonconstant
meromorphic functions
$j\colon \Cp1\to\Cp1$ satisfying
conditions \loccit{degree}--\loccit{mod2} above results in a
fiberwise deformation of the corresponding
\rom(\via\ Proposition~\ref{j->B}\rom)
simplified trigonal curves.
\endproposition

\remark
Condition~\iref{j.deformations}{mod3} means that a root of~$j$ of
multiplicity divisible by~$3$ may join another root and,
conversely, a root of large multiplicity may break into several
roots, all but one having multiplicities divisible by~$3$.
Condition~\iref{j.deformations}{mod2} should be interpreted similarly.
\endremark

\remark
Note that just an equisingular (not necessarily fiberwise)
deformation of trigonal curves does not always result in a
deformation of their $j$-invariants. In the case of simplified
curves, the degree of~$j_B$ drops whenever
a type~$\tilde\bA_0^*$ singular fiber of~$B$ joins another
singular fiber, of type~$\tilde\bA_0^*$, $\tilde\bA_1$,
or~$\tilde\bA_2$, to form a fiber of type~$\tilde\bA_0^{**}$,
$\tilde\bA_1^*$, or~$\tilde\bA_2^*$, respectively,
see Table~\ref{table.j}.
\endremark

\proof
The direct statement follows essentially from Table~\ref{table.j}.
Indeed, the multiplicities of the poles of~$j_B$,
$(\omod3)$-multiplicities of its roots, and
$(\omod2)$-multiplicities of the roots of $j_B-1$ are encoded in
the singular fibers of~$B$, and the degree $\deg j_B$ can be found
as the sum of the multiplicities of all poles of~$j_B$. Since
the expression for~$j_B$ depends `continuously' on the
coefficients of the Weierstrass equation and $\deg j_B$ remains
constant, there is no extra cancellation during the deformation
and the map $j_B\colon \Cp1\to\Cp1$ changes continuously.

The converse statement follows from the construction of the
simplified trigonal curve~$B$ from a given $j$-invariant~$j$, see
the proof of Proposition~\ref{j->B}. Since the degree
$\deg l=\deg j$ remains constant, the polynomials~$p$ and~$q$ in
the irreducible representation $l=p/q$ change continuously during
the deformation. Crucial is the fact that the passage from $p/q$
to $\bar p/\bar q$ depends only on the roots of~$p$ and~$q$ whose
multiplicity is not divisible by~$3$ and~$2$, respectively. Hence,
due to Conditions~\iref{j.deformations}{mod3} and~\ditto{mod2},
the degree $\deg\bar p=\deg\bar q$ will remain constant, the
polynomials~$\bar p$ and~$\bar q$ will change continuously, and so
will the coefficients~$g_2$, $g_3$ of the Weierstrass equation.
The fact that the resulting deformation of the trigonal curves is
fiberwise follows again from Table~\ref{table.j}.
\endproof

\section{Dessins d'enfants and skeletons}\label{S.dessins}

According to Propositions~\ref{j->B} and~\ref{j.deformations}, the
study of simplified trigonal curves in Hirzebruch surfaces is
reduced to the study of meromorphic functions $j\colon \Cp1\to\Cp1$
with three `essential' critical values~$0$, $1$, and~$\infty$ and,
possibly, a few other critical values. Following
S\period~Orevkov~\cite{Orevkov}, \cite{Orevkov-private}, we
employ a modified version of Grothendieck's \emph{dessins
d'enfants}. Below, we outline briefly the basic concepts and
principal results; for more details and proofs we refer
to~\cite{DIK.elliptic}, Sections~5.1 and~5.2.
Note that \cite{DIK.elliptic} deals with a
real version of the theory, where functions (graphs) are supplied
with an anti-holomorphic (respectively, orientation reversing)
involution; however, all proofs apply to the settings of this
paper literally, with the real structure ignored.

Since, in this paper, we deal with rational ruled surfaces only,
we restrict the further exhibition to the case of graphs in the
sphere $\base\cong\Cp1$.

\subsection{Trichotomic graphs}\label{trichotomic}
Given a graph~$\Gamma\subset\base$,
we denote by~$\Base\Gamma$ the closed cut
of~$\base$ along~$\Gamma$. The connected components of~$\Base\Gamma$
are called the \emph{regions} of~$\Gamma$.
(Unless specified otherwise, in the topological
part of this section we are working in the \PL-category.)

A \emph{trichotomic graph}
is an embedded oriented
graph $\Gamma\subset\base$
decorated with the
following additional structures (referred to as \emph{colorings}
of the edges and vertices of~$\Gamma$, respectively):
\roster
\item"--"
each edge of~$\Gamma$ is of one of the three kinds: \solid-, \bold-,
or \dotted-;
\item"--"
each vertex of~$\Gamma$ is of one of the four kinds: \black-, \white-,
\cross-, or monochrome (the vertices of the first three kinds being called
\emph{essential})
\endroster
and satisfying the following conditions:
\roster
\item\local{tg-valency}
the valency of each essential vertex of~$\Gamma$ is at
least~$2$,
and the valency of each monochrome vertex of~$\Gamma$ is at
least~$3$;
\item\local{tg-oriented}
the orientations of the edges of~$\Gamma$ form an orientation
of the boundary $\partial\Base\Gamma$;
this orientation extends to an orientation of~$\Base\Gamma$;
\item\local{tg-monochrome}
all edges incident to a monochrome vertex are of the same kind;
\item\local{tg-cross}
\cross-vertices are incident to incoming \dotted-edges and
outgoing \solid-edges;
\item\local{tg-black}
\black-vertices are incident to incoming \solid-edges and
outgoing \bold-edges;
\item\local{tg-white}
\white-vertices are incident to incoming \bold-edges and
outgoing \dotted-edges.
\endroster
In \loccit{tg-cross}--\loccit{tg-white} the lists are complete,
\ie, vertices cannot be incident to edges of other kinds or with
different orientation.

Condition~\loccit{tg-oriented} implies that the
orientations of the edges incident to a vertex alternate.
In particular, all vertices of~$\Gamma$
have even valencies.

\subsection{Dessins}\label{dessins}
In view of~\iref{trichotomic}{tg-monochrome}, the monochrome vertices of a
trichotomic graph $\Gamma$ can further be subdivided into \solid-,
\bold-, and \dotted-, according to their incident edges. A path
in~$\Gamma$ is called \emph{monochrome} if all its vertices are
monochrome. (Then, all vertices of the path are of the same kind,
and all its edges are of the same kind as its vertices.) Given two
monochrome vertices $u,v\in\Gamma$, we say that $u\prec v$ if
there is an oriented monochrome path from~$u$ to~$v$. (Clearly,
only vertices of the same kind can be compatible.) The graph is
called \emph{admissible} if $\prec$ is a partial order. Since
$\prec$ is obviously transitive, this condition is equivalent to
the requirement that $\Gamma$ should have no oriented monochrome
cycles.

In this paper, an admissible trichotomic graph is
called a \emph{dessin}.

\remark
Note
that the orientation of~$\Gamma$ is almost
superfluous. Indeed, $\Gamma$ may have at most two orientations
satisfying~\iref{trichotomic}{tg-oriented},
and if $\Gamma$ has at least one
essential vertex, its orientation is uniquely determined by
\iref{trichotomic}{tg-cross}--\ditto{tg-white}.
Note also that each connected component of an admissible
graph does have essential vertices (of all three kinds),
as otherwise any component
of~$\partial\Base\Gamma$ would be an oriented monochrome cycle.
\endremark

\remark
In fact, all three decorations of a dessin~$\Gamma$
(orientation and the two colorings) can be recovered from any of
the colorings. However, for clarity we retain both
colorings in the diagrams.
\endremark

\subsection{The dessin of a trigonal curve}\label{j->Gamma}
Any orientation preserving
ramified covering $j\colon \base\to\Cp1$ defines a trichotomic
graph $\Gamma(j)\subset\base$. As a set, $\Gamma(j)$ is the
pull-back $j^{-1}(\Rp1)$. (Here, $\Rp1\subset\Cp1$ is the fixed
point set of the standard real structure $z\mapsto\bar z$.)
The trichotomic graph structure
on~$\Gamma(j)$ is introduced as follows: the \black--, \white--,
and \cross-vertices are  the pull-backs of~$0$, $1$, and~$\infty$,
respectively (monochrome vertices being the ramification points
with other real critical values), the edges are \solid-, \bold-,
or \dotted- provided that their images belong to $[\infty,0]$,
$[0,1]$, or $[1,\infty]$, respectively, and the orientation
of~$\Gamma(j)$ is that induced from the positive orientation
of~$\Rp1$ (\ie, order of~$\R$).

As shown in~\cite{DIK.elliptic}, a trichotomic graph
$\Gamma\subset\base$ is a dessin if and only if it has the form
$\Gamma(j)$ for some orientation preserving
ramified covering $j\colon \base\to\Cp1$; the latter is determined
by~$\Gamma$ uniquely
up to homotopy in the class of ramified coverings
having a fixed trichotomic graph.

We define the \emph{dessin} $\Gamma(B)$ of a trigonal curve~$B$
as the dessin
$\Gamma(j)$ of its $j$-invariant $j\colon \Cp1\to\Cp1$. The
correspondence between the singular fibers of a simplified
trigonal curve~$B$ and the vertices of its dessin $\Gamma(B)$ is
given by Table~\ref{table.j}
(see $j(z)$\,), the valency of a vertex~$z$ being
twice the ramification index $\ind_zj$.
The \black-- (respectively, \white--) vertices of
$\Gamma(B)$ of valency $0\bmod6$
(respectively, $0\bmod4$) correspond to the nonsingular fibers
of~$B$ with complex multiplication of order~$3$
(respectively,~$2$); such vertices are called \emph{nonsingular},
whereas all other essential vertices of~$\Gamma$ are called
\emph{singular}.

\subsection{Equivalence of dessins}\label{dessin.equivalence}
Let $\Gamma\subset\base$ be a trichotomic graph, and let~$v$ be a
vertex of~$\Gamma$. Pick a regular neighborhood $U\subset\base$
of~$v$ and replace the intersection $\Gamma\cap U$ with another
decorated graph, so that the result~$\Gamma'$ is again a
trichotomic graph. If $\Gamma'\cap U$ contains essential
vertices of at most one kind, then $\Gamma'$ is called a
\emph{perturbation} of~$\Gamma$ (at~$v$),
and the original graph~$\Gamma$ is called a
\emph{degeneration} of~$\Gamma'$.

A perturbation~$\Gamma'$ of a dessin is also a dessin if and only if
the intersection $\Gamma'\cap U$ contains no oriented monochrome cycles.
There
are
no simple local criteria
for the admissibility of a degeneration.

\remark
Assume that the perturbation~$\Gamma'$ is a dessin.
Since the intersection $\Gamma'\cap\partial U$ is
fixed, the assumption on $\Gamma'\cap U$ implies that
$\Gamma'\cap U$ either is monochrome (if $v$ is monochrome) or
consists of monochrome vertices, essential vertices of the same
kind as~$v$, and edges of the two kinds incident to~$v$.
\endremark

A perturbation~$\Gamma'$ of a dessin~$\Gamma$ at a vertex~$v$ (and
the inverse degeneration of~$\Gamma'$ to~$\Gamma$) is called
\emph{equisingular} if $v$ is not a \cross-vertex and the
intersection $\Gamma'\cap U$ contains at most one singular
\black-- or \white-vertex.
Two dessins $\Gamma',\Gamma''\subset\base$
are said to be \emph{equivalent} if they can be connected by a
chain $\Gamma'=\Gamma_0,\Gamma_1,\ldots,\Gamma_n=\Gamma''$ of
dessins, where each~$\Gamma_i$, $1\le i\le n$, either is isotopic
to~$\Gamma_{i-1}$
or is an equisingular perturbation or
degeneration of~$\Gamma_{i-1}$. Clearly, equivalence of dessins is
an equivalence relation.

\remark
By an \emph{isotopy} between two dessins~$\Gamma'$ and~$\Gamma''$ we
mean a \PL-family $\phi_t$
of \PL-autohomeomorphisms of~$\base$
such that $\phi_0=\id$ and $\phi_1(\Gamma')=\Gamma''$, the latter
map taking vertices to vertices and edges to edges and preserving
both colorings of the dessins. Note that, since the mapping class
group of~$\base$ is trivial, one can just require that $\Gamma'$
is taken to~$\Gamma''$ by an \emph{orientation preserving}
\PL-autohomeomorphism of~$\base$ respecting the graph structure
and the colorings.
\endremark

The following statement,
essentially based on the Riemann existence theorem,
is an immediate consequence of
Propositions~\ref{j->B} and~\ref{j.deformations} and the results
of~\cite{DIK.elliptic} (particularly, Corollaries~5.1.8 and~5.2.3,
with the real structure ignored).

\theorem\label{B->Gamma}
The map $B\mapsto\Gamma(B)$
sending a trigonal curve~$B$ to its dessin
establishes a one-to-one correspondence between the set of
fiberwise deformation classes of simplified
trigonal curves in Hirzebruch
surfaces and the set of equivalence classes of dessins.
\qed
\endtheorem

\definition[Maximal curves and dessins]\label{def.maximal}
A dessin $\Gamma\subset\base$ is called \emph{maximal} if it
satisfies the following conditions:
\roster
\item\local{mono}
all vertices of~$\Gamma$ are essential;
\item\local{valency}
all \black-- (respectively, \white--) vertices of~$\Gamma$ have
valency~$\le6$ (respectively,~$\le4$);
\item\local{triangles}
all regions of~$\Gamma$ are triangles.
\endroster
A simplified trigonal curve~$B$ is called \emph{maximal}
if its dessin $\Gamma(B)$ is maximal.
\enddefinition

\remark
Conditions~\iref{def.maximal}{mono} and~\ditto{triangles} in the
definition of a maximal dessin
can be restated as the requirement that
the function $j\colon \base\to\Cp1$ constructed from~$\Gamma$,
see~\ref{j->Gamma},
should have no critical values other than $0$, $1$, and~$\infty$.
\endremark

\remark
Any maximal trigonal curve is defined over an algebraic number
field. Indeed, as any function with three critical values, the
rational function $j\colon \Cp1\to\Cp1$ has finitely many Galois
conjugates and hence
is defined over an algebraic
number field. Then, the construction in the proof of
Proposition~\ref{j->B} shows that the coefficients~$g_2$, $g_3$ of
the Weierstra{\ss} equation are defined over the splitting field
of~$j$. (One may need to add to the field some roots and poles
of~$j$.)
\endremark

\proposition\label{max.B->Gamma}
A trigonal curve~$B'$ is fiberwise deformation
equivalent to a maximal trigonal curve~$B$ if and only if the
dessins $\Gamma(B')$ and $\Gamma(B)$ are isotopic.
Furthermore, a permutation of the singular fibers of a maximal
trigonal curve~$B$ is realized by a fiberwise self-deformation if and
only if the corresponding permutation of the vertices of
$\Gamma(B)$
is induced by an isotopy of~$\Gamma(B)$.
\endproposition

\proof
A maximal dessin does not admit nontrivial equisingular
perturbations (due to Conditions~\loccit{mono}
and~\loccit{valency}, as an equisingular perturbation requires
a vertex of high valency) or degenerations (due to
Condition~\loccit{triangles}, as a perturbation produces more
than triangle regions). Hence, any equivalence to a maximal dessin
is an isotopy.
\endproof

\proposition\label{no.degeneration}
A trigonal curve~$B$ is maximal if and only if it does not admit
a nontrivial degeneration,
see~\ref{deformations}, to a non-isotrivial curve.
\endproposition

\proof
Let $\Gamma=\Gamma(B)$. After a small deformation, we can assume
that $\Gamma$ satisfies the general position
assumptions~\iref{def.maximal}{mono} and~\ditto{valency}. Then,
if~$B$ is not maximal, $\Gamma$ has a region~$R$ whose boundary
contains at least two \cross-vertices, and these two vertices can
be brought together within~$R$. This degeneration of~$\Gamma$
results in a nontrivial degeneration of the curve.

Conversely, assume that $B$ has a nontrivial degeneration to a
curve~$B_0$, which is necessarily trigonal. Then, up to isotopy,
$\Gamma$ is obtained from $\Gamma(B_0)$ by removing disjoint
regular neighborhoods of some of its vertices and replacing them
with new decorated graphs. (Since $\deg j$ may change, it is no
longer required
that each of the new graphs should contain essential vertices of
at most one kind. Note that we do not discuss the realizability of
any such modification by an actual degeneration of curves.)
If this procedure is nontrivial, it
results in a graph~$\Gamma$ with at least
one non-triangular region.
\endproof

\subsection{Skeletons}\label{skeletons}
An \emph{abstract skeleton} is a connected
planar map $\Sk\subset\base$ whose
vertices have valencies at most three; we allow the possibility of
\emph{hanging edges}, \ie, edges with only one end attached to a
vertex. An \emph{isomorphism} between two abstract skeletons~$\Sk'$
and~$\Sk''$ is an orientation preserving \PL-autohomeomorphism
of~$\base$ taking~$\Sk'$ to~$\Sk''$.

The \emph{skeleton} of a maximal trigonal curve~$B$ is the
skeleton $\Sk(B)\subset\base$
obtained from the dessin $\Gamma(B)$ by removing
all \cross-vertices and incident edges (\ie, all \solid- and \dotted-
edges) and disregarding the \white-vertices. (Note that the
resulting graph is indeed connected due to
Condition~\iref{def.maximal}{triangles}.)
Clearly, a maximal
dessin~$\Gamma$ is uniquely (up to homotopy) recovered from its
skeleton~$\Sk$: one should place a \white-vertex at the middle of
each edge (at the free end of each hanging edge), place a
\cross-vertex $v_R$ at the center of each region~$R$ of~$\Sk$,
and connect this vertex
$v_R$ to the \black-- and \white-vertices in the
boundary~$\partial R$ by appropriate (respectively, \solid- and
\dotted-) edges. The last operation is unambiguous as, due to the
connectedness of~$\Sk$, each open region~$R$ is a topological disk.

\theorem\label{B->Sk}
The map $B\mapsto\Sk(B)$
sending a maximal trigonal curve~$B$ to its skeleton
establishes a one-to-one correspondence between the set of
fiberwise deformation classes of maximal
trigonal curves in Hirzebruch
surfaces and the set of isomorphism classes of abstract skeletons
in~$\base$.
\endtheorem

\proof
The statement follows from the
correspondence between maximal dessins and skeletons described above,
Theorem~\ref{B->Gamma}, and Proposition~\ref{max.B->Gamma}.
\endproof

\remark
Removing from a dessin~$\Gamma(j)$ all \cross-vertices and incident
edges results in a classical \emph{dessin d'enfants} in the
sense of Grothendieck, \ie, the bipartite graph obtained as the
pull-back $j^{-1}([0,1])$. The passage to the skeletons is a
further simplification due to the fact that, under the
assumptions on maximal dessins,
all \white-vertices have valency at most two.
\endremark

\remark
Theorem~\ref{B->Sk} suggests that, in general, the classification
of maximal trigonal curves with a prescribed combinatorial type of
singular fibers is a wild problem: one would have to enumerate all
planar maps with prescribed valencies of vertices and numbers of
edges of regions.
The only general result in this direction that I am aware of is the
Hurwitz formula~\cite{Hurwitz} (see also~\cite{Bauer}), which
establishes a relation between a certain weighed count of planar
maps (more precisely, ramified coverings of~$\Cp1$, not necessarily
connected)
and characters of symmetric groups.
\endremark

\subsection{Vertex count}\label{counts}
We conclude this section with a few simple counts. For a
dessin~$\Gamma$, denote by $\nstar=\nstar(\Gamma)$ the total number of
\star-vertices (where \star-\ is either \black-, or \white-, or
\cross-), and by $\nstar(i)$, $i\in\N$,
the number of \star-vertices of
valency~$2i$. (Recall that valencies of all vertices of a dessin
are even.) Consider a trigonal curve $B\subset\Sigma_k$, its
$j$-invariant $j\colon \Cp1\to\Cp1$, and its dessin $\Gamma=\Gamma(B)$.
Counting the number of points in one of the three special fibers
of~$j$, one obtains
\begin{equation}
\deg j=\sum_{i>0}i\nblack(i)
 =\sum_{i>0}i\nwhite(i)
 =\sum_{i>0}i\ncross(i).\label{eq.deg.j}
\end{equation}
Since $B$ can be perturbed to a generic trigonal curve in the same
surface~$\Sigma_k$, and a generic curve has $\deg\Delta=6k$
simplest singular fibers, Table~\ref{table.j} yields
\begin{equation}
6k=\sum_{i>0}i\ncross(i)
 +2\sum_{i=1(3)}\nblack(i)
 +3\sum_{i=1(2)}\nwhite(i)
 +4\sum_{i=2(3)}\nblack(i).\label{eq.6k}
\end{equation}
(Alternatively, one can notice that the first term
in~\eqref{eq.6k} equals $\deg j$, see~\eqref{eq.deg.j}, and the
remaining part of the sum is $6k-\deg j$, see $\Delta\deg j$ in
Table~\ref{table.j}.) Finally, the Riemann--Hurwitz formula
applied to~$j$ results in the inequality
\begin{equation}
\nblack+\nwhite+\ncross\ge\deg j+2,\label{eq.R-H}
\end{equation}
which turns into an equality if and only if $j$ has no critical
values other than~$0$, $1$, and~$\infty$, \ie,
Conditions~\iref{def.maximal}{mono} and~\ditto{triangles} are
satisfied.

\section{The fundamental group}\label{S.group}

\subsection{The braid group}\label{B3}
Recall that the \emph{braid group} $\BG3$ can be defined as the
group of automorphisms of the free group $G=\<\Ga_1,\Ga_2,\Ga_3>$
sending each generator to a conjugate of another generator and
leaving the product $\Ga_1\Ga_2\Ga_3$ fixed. We assume that the
action of~$\BG3$ on~$G$ is \emph{from the left}.
One has
$\BG3=\<\Gs_1,\Gs_2\,|\,\Gs_1\Gs_2\Gs_1=\Gs_2\Gs_1\Gs_2>$, where
$$
\Gs_1\colon (\Ga_1,\Ga_2,\Ga_3)\mapsto(\Ga_1\Ga_2\Ga_1^{-1},\Ga_1,\Ga_3),
\quad
\Gs_2\colon (\Ga_1,\Ga_2,\Ga_3)\mapsto(\Ga_1,\Ga_2\Ga_3\Ga_2^{-1},\Ga_2).
$$
We will also consider the elements $\Gs_3=\Gs_1^{-1}\Gs_2\Gs_1$
and
$\tau=\Gs_2\Gs_1=\Gs_3\Gs_2=\Gs_1\Gs_3$. The center of~$\BG3$ is
the infinite cyclic group generated by~$\tau^3$.

Note that the maps
$(\Gs_1,\Gs_2)\mapsto(\Gs_2,\Gs_3)\mapsto(\Gs_3,\Gs_1)$ define
automorphisms of~$\BG3$; in particular, the pairs $(\Gs_2,\Gs_3)$
and $(\Gs_3,\Gs_1)$ are subject to all relations that hold for
$(\Gs_1,\Gs_2)$. In what follows, we use the convention
$\Gs_{3l+i}=\Gs_i$, $i=1,2,3$, $l\in\Z$.

The \emph{degree} $\deg\Gb$ of a braid $\Gb\in\BG3$ is defined as
its image under the abelinization homomorphism $\BG3\to\Z$,
$\Gs_1,\Gs_2\mapsto1$. A braid is uniquely recovered from its
degree and its image in the quotient $\BG3/\tau^3$.

\subsection{Van Kampen's method}\label{vanKampen.method}
Let $B\subset\Sigma=\Sigma_k$ be a generalized trigonal curve, and let
$E\subset\Sigma$ be the exceptional section.
The fundamental group $\pi_1(\Sigma\sminus(B\cup E))$ can be found
using an analogue of van Kampen's method~\cite{vanKampen}
applied to the ruling of~$\Sigma$.
Pick a fiber~$F_\infty$ (singular or not) over a point
$\infty\in\Cp1$
and trivialize the ruling over
$\Cp1\sminus\infty$.
Let
$F_1,\ldots,F_r$ be the singular fibers of~$B$ other
than~$F_\infty$.
Pick a nonsingular fiber~$F$ distinct from~$F_\infty$ and
a generic section~$S$
disjoint from~$E$ and intersecting all fibers
$F,F_1,\ldots,F_r,F_\infty$ outside of~$B$.

Clearly,
$F\sminus(B\cup E)$ is the plane $\C^1=F\sminus E$ with three punctures.
Consider the group
$G=\pi_1(F\sminus(B\cup E),F\cap S)$, and let $\Ga_1$,
$\Ga_2$, $\Ga_3$ be a standard set of generators of~$G$. Let,
further, $\Gg_1,\ldots\Gg_r$ be a standard set of generators of the
fundamental group
$\pi_1(S\sminus(F_\infty\cup\bigcup_{j=1}^rF_r),S\cap F)$,
so that $\Gg_j$ is a
loop around~$F_j$, $j=1,\ldots,r$. For each $j=1,\ldots,r$,
dragging the fiber~$F$ along~$\Gg_j$ and keeping the base
point in~$S$ results in a certain automorphism $m_j\colon G\to G$,
called the \emph{braid monodromy} along~$\Gg_j$. Strictly
speaking, $m_j$ is not necessarily a braid (unless $B$ is
disjoint from~$E$); however, it still has the property that the
image $m_j(\Ga_i)$ of each standard generator~$\Ga_i$, $i=1,2,3$,
is a conjugate of another generator~$\Ga_{i'}$.

According to van Kampen,
the group
$\pi_1(\Sigma\sminus(B\cup E\cup F_\infty\cup\bigcup_{j=1}^rF_r),S\cap F)$
is given by the representation
$$
\bigl<\Ga_1,\Ga_2,\Ga_3,\Gg_1,\ldots,\Gg_r\bigm|
 \text{$\Gg_j^{-1}\Ga_i\Gg_j=m_j(\Ga_i)$, $i=1,2,3$,
 $j=1,\ldots,r$}\bigr>,
$$
and patching back a fiber~$F_j$, $j=1,\ldots,r$, results in an
additional relation $\Gg_j=1$. Thus,
if $B$ has no linear components,
the resulting representation
for the group
$\pi_1(\Sigma\sminus(B\cup E\cup F_\infty),S\cap F)$ is
$$
\bigl<\Ga_1,\Ga_2,\Ga_3\bigm|
 \text{$\Ga_i=m_j(\Ga_i)$, $i=1,2,3$, $j=1,\ldots,r$}\bigr>.
$$
Patching back the remaining fiber~$F_\infty$ gives one more relation
$\Gg=1$, where $\Gg$ is the class of a small loop in~$S$ around
$S\cap F_\infty$; an expression of~$\Gg$ in terms of $\Ga_1$, $\Ga_2$,
$\Ga_3$ in the special case of trigonal curves is found below,
see Remark in~\ref{proper.section}.

\remark
Van Kampen's approach applies as well in the case when the curve
has linear components:
for each such component,
one should keep the corresponding generator~$\Gg_j$ and keep the
relation $\Gg_j^{-1}\Ga_i\Gg_j=m_j(\Ga_i)$ instead of
$\Ga_i=m_j(\Ga_i)$.
\endremark

\proposition\label{pi1}
Let $B_k\subset\Sigma_k$ and $B_{k+1}\subset\Sigma_{k+1}$ be two
generalized trigonal curves, so that $B_k$ is obtained
from~$B_{k+1}$ by an elementary transformation whose blow-up
center~$O$ does not belong to~$B_{k+1}$. Then there is a natural
isomorphism
$$
\pi_1(\Sigma_k\sminus(B_k\cup E_k))=
 \pi_1(\Sigma_{k+1}\sminus(B_{k+1}\cup E_{k+1})),
$$
where $E_i\subset\Sigma_i$, $i=k,k+1$, are the exceptional sections.
\endproposition

\proof
Let $F_{k+1}\subset\Sigma_{k+1}$ be the fiber through~$O$, and let
$F_k\subset\Sigma_k$ be the fiber contracted by the
inverse elementary transformation. The
diffeomorphism~\eqref{elementary.diffeo} induces an isomorphism
$$
\pi_1(\Sigma_k\sminus(B_k\cup E_k\cup F_k))=
 \pi_1(\Sigma_{k+1}\sminus(B_{k+1}\cup E_{k+1}\cup F_{k+1})).
$$
The group $\pi_1(\Sigma_k\sminus(B_k\cup E_k))$ is obtained from
$\pi_1(\Sigma_k\sminus(B_k\cup E_k\cup F_k))$ by adding the
relation $[\partial\Gamma_k]=1$, where $\Gamma_k\subset\Sigma_k$
is a small analytic disk transversal to~$F_k$ and disjoint from
all other curves involved. Similarly, patching the fiber~$F_{k+1}$ results
in an additional relation $[\partial\Gamma_{k+1}]=1$, where
$\Gamma_{k+1}\subset\Sigma_{k+1}$ is a small analytic disk
transversal to~$F_{k+1}$ and disjoint from the other curves
in~$\Sigma_{k+1}$. Under the assumptions, one can
choose~$\Gamma_{k+1}$ passing through the blow-up center~$O$; then
its proper transform can be taken for~$\Gamma_k$. Hence, one has
$[\partial\Gamma_k]=[\partial\Gamma_{k+1}]$, and the two quotient
groups are isomorphic.
\endproof

\proposition\label{abelian.group}
Let $C\subset\Cp2$ be an algebraic curve of degree~$m$ with
a distinguished singular point~$O$ of
multiplicity~$(m-3)$ and without linear components.
Assume that $C$ has a branch~$b$ at~$O$ of
type~$\bE_{12}$. Then $C$ is irreducible and
the fundamental group $\pi_1(\Cp2\sminus C)=\Z_m$ is abelian.
\endproposition

\proof
Blow~$O$ up and consider the proper transform
$B_1\subset\Sigma_1$ of~$C$,
see~\ref{trigonal.model}. The transform of~$b$ is a type~$\bE_6$
singular point of~$B_1$, and the elementary transformation centered
at this point converts~$B_1$ to a generalized trigonal curve
$B_2\subset\Sigma_2$ with a type~$\tilde\bA_0^{**}$ singular
fiber. In particular, the curve is irreducible.

The inverse transformation is as in Proposition~\ref{pi1},
\ie, its blow-up center does not belong to the curve~$B_2$ or the
exceptional section~$E_2$. Hence, one has
$$
\pi_1(\Cp2\sminus C)=\pi_1(\Sigma_1\sminus(B_1\cup E_1))=
 \pi_1(\Sigma_2\sminus(B_2\cup E_2)).
$$
(The first
isomorphism is obvious; the second one is given by
Proposition~\ref{pi1}.) The last group can be found using van
Kampen's method, see~\ref{vanKampen.method}.
Under an appropriate choice of the
generators $\Ga_1$, $\Ga_2$, $\Ga_3$, the braid monodromy~$m$
about a type~$\tilde\bA_0^{**}$ singular fiber is $\tau\in\BG3$,
and the relations $m(\Ga_i)=\Ga_i$, $i=1,2,3$, yield
$\Ga_1=\Ga_2=\Ga_3$. Hence, the group is abelian.
\endproof

\proposition\label{abelian.group.ex}
Let $C$ be the union of an irreducible curve as in
Proposition~\ref{abelian.group} and $r\ge1$ linear components none
of which is tangent to the branch~$b$ of type~$\bE_{12}$. Then one
has
$\pi_1(\Cp2\sminus C)=\Z\times\<\Gg_1,\ldots\Gg_{r-1}>$. In
particular, if $r\le2$, the group is still abelian.
\endproposition

\proof
As in the proof of Proposition~\ref{abelian.group},
there is a relation $\Ga_1=\Ga_2=\Ga_3$, and due to the
properties of the
braid monodromy (each generator is taken to a conjugate of a
generator) the relations
$\Gg_j^{-1}\Ga_i\Gg_j=m_j(\Ga_i)$ turn into $[\Gg_j,\Ga_i]=1$.
\endproof

\section{The braid monodromy}\label{S.monodromy}

In this section, we describe the braid monodromy of a
simplified trigonal curve.
We fix such a curve $B\subset\Sigma=\Sigma_k$ and let
$\Gamma=\Gamma(B)$. Further, we
denote by~$F_z$ the fiber over a point $z\in\Cp1$, and let
$B_z=B\cap F_z$ and $E_z=E\cap F_z$, where $E\subset\Sigma$ is the
exceptional section.
Note that
$F_z\sminus E_z$ is an affine space over~$\C^1$; in particular, one
can speak about its orientation,
lines, circles, angles, and length ratios. We use the notation $F_z^\circ$
for the punctured plane $F_z\sminus(B_z\cup E_z)$.

\subsection{Geometry of the fibers}\label{tg->topology}
The definition of the $j$-invariant gives an easy way to
recover the topology of~$B$ from its dessin $\Gamma$.
The set~$B_z$ consists of a single triple point if $z$ is a
singular \black-- or \white-vertex. If $z$ is
a \cross-vertex, $B_z$ consists of two points, one simple and one
double. In all other cases, $B_z$ consists of three simple points,
whose position in~$F_z\sminus E_z$ can be characterized as
follows.
\roster
\item\local{b.region}
If $z$ is an inner point of a region of~$\Gamma$, the three points
of~$B_z$ form a triangle with all three
edges distinct. Hence, the restriction of the projection
$B\to\Cp1$ to the
interior of each region of~$\Gamma$ is a trivial covering.
\item\local{b.dotted}
If $z$ belongs to a \dotted- edge of~$\Gamma$, the three points
of~$B_z$ are collinear. The ratio
$(\text{smallest distance})/(\text{largest distance})$
is in $(0,\frac12)$; it
tends to~$0$ (respectively, $\frac12$) when $z$ approaches a \cross--
(respectively,~\white--) vertex.
\item\local{b.edge}
If $z$ belongs to a \solid- (\bold-) edge of~$\Gamma$, the three
points of~$B_z$ form an isosceles triangle
with the angle at the vertex less than (respectively, greater than)
$\pi\!/3$. The angle
tends to~$0$, $\pi\!/3$, or~$\pi$ when $z$ approaches, respectively,
a \cross--, \black--, or \white-vertex.
\endroster
Furthermore, a simple model example proves the following statement.
\roster
\item[4]\local{b.orientation}
For a point~$z$ as in~\loccit{b.region}, arrange the vertices
of~$B_z$ in ascending order based on the length of
the opposite edge.
The resulting orientation of~$B_z$ is
counterclockwise if and only if
$\Im j_B(z)>0$.
\endroster

\subsection{Proper sections}\label{proper.section}
To define the braid monodromy, we need to fix a `fiber at
infinity'~$F_\infty$, see~\ref{vanKampen.method}, and
a generic section~$S$
that would provide the base points
$S_z=S\cap F_z\in F_z^\circ$.
We take for~$F_\infty$ the fiber over a fixed point
$\infty\notin\Gamma$, and construct~$S$ as a small perturbation of
$E+kF'$, where $F'$ is the fiber over a point~$z'$ in the
same open region of~$\Gamma$ as~$\infty$. If the perturbation is
sufficiently small, the section~$S$ has the following property:
{\proclaimfont there is a closed neighborhood $K\ni\infty$ disjoint
from~$\Gamma$ and such that, for each point
$z\in\Cp1\sminus K$, the base point $S_z\in F_z$ is outside a
disk $U_z\subset F_z$ containing~$B_z$ and centered at its
barycenter \rom(\cf.~Figure~\ref{fig.basis}, right, below\rom).}
In what
follows, a section~$S$ satisfying this property is called
\emph{proper} and, when speaking about the fundamental group
$\pi_1(F_z^\circ,S_z)$, we always assume that the
point~$z$ is outside the above closed neighborhood~$K$.

Note that, together with the exceptional section~$E$ and the zero
section given by $z\mapsto\text{(the barycenter of~$B_z$)}$, a
proper section~$S$ gives a trivialization of the ruling over
$\Cp1\sminus K$, which is necessary to define the braid monodromy.

\remark
From the construction of a proper section~$S$ it follows that the
class~$\Gg$ of a small loop in~$S$ surrounding~$F_\infty\cap S$
(see~\ref{vanKampen.method}) is, up to conjugation, given by
$\Gg=(\Ga_1\Ga_2\Ga_3)^k\Gg_1\ldots\Gg_r$. Hence, in this case,
the final relation in van Kampen's method is
$(\Ga_1\Ga_2\Ga_3)^k=1$.
\endremark

\subsection{Markings and canonical bases}\label{basis}
Let $z\in\Gamma$ be a nonsingular \black-vertex.
According to~\iref{tg->topology}{b.edge}, the three points
of the set~$B_z$ form an equilateral triangle. There is a
natural one-to-one correspondence between the \bold- edges
incident to~$z$ and the points of~$B_z$: an edge~$e$ corresponds
to the point $p\in B_z$ that turns into the vertex of the isosceles
triangle when $z$ slides from its original position along~$e$. In
fact, the same point~$p$ turns into the vertex of the isosceles
triangle when $z$ slides along the \solid- edge~$e'$ opposite
to~$e$, so that the two other points are brought together over the
\cross-vertex ending~$e'$.

\begin{figure}[here]
\centerline{\epsffile{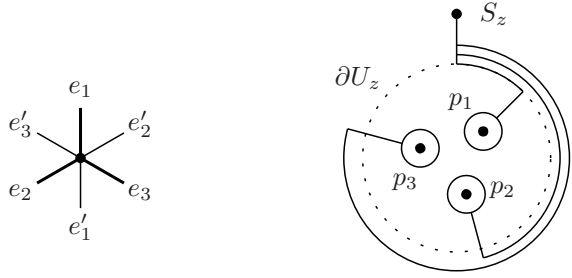}}
\caption{A canonical basis for~$G_z$}\label{fig.basis}
\end{figure}

In what follows, we always assume that the three bold edges
$e_1,e_2,e_3$ incident to~$z$ are oriented
in the \emph{counterclockwise} direction,
as in Figure~\ref{fig.basis}, left. Such an ordering is called a
\emph{marking} at~$z$, and an edge~$e_i$ incident to~$z$ is said
to have \emph{index}~$i$ at~$z$. A marking at~$z$ is uniquely
determined by assigning an index to one of the three \bold- edges
incident to~$z$. Alternatively, a marking is determined by
assigning an index to one of the three points constituting~$B_z$.

A \emph{marking} of a dessin~$\Gamma$ is defined as a collection
of markings at each nonsingular \black-vertex of~$\Gamma$. The
notion of marking and index of edges extends to skeletons in the
obvious way.

Using~\iref{tg->topology}{b.region}--\ditto{b.edge},
from~\iref{tg->topology}{b.orientation} it follows that
{\proclaimfont if $e_1, e_2, e_3$ is a marking at a nonsingular
\black-vertex~$z$, the
corresponding points $p_1,p_2,p_3\in B_z$ form the clockwise
orientation of the triangle~$B_z$ \rom(Figure~\ref{fig.basis},
right\rom).}

Pick a proper section~$S$, see~\ref{proper.section}, and consider
the group $G_z=\pi_1(F_z^\circ,S_z)$.
A \emph{canonical basis} for~$G_z$ is a basis
$\Ga_1,\Ga_2,\Ga_3$ shown in Figure~\ref{fig.basis}, right,
where the space
$F_z^\circ$ is regarded as the affine line
$F_z\sminus E_z$ punctured at $p_1,p_2,p_3\in B_z$.
More precisely, each
element~$\Ga_i$ is the class of the loop formed by a small
counterclockwise circle about~$p_i$, $i=1,2,3$, which is connected
to~$S_z$ by a radial segment, an arc of a circle $\partial U_z$
separating~$S_z$ from~$B_z$ (\cf.~\ref{proper.section}), and
another radial segment, common for all three loops. It is required
that each consecutive arc is $2\pi\!/3$ longer than the previous
one; however, we do not make any assumption about the length of
the first arc: it is defined up to a multiple of~$2\pi$. As a
result, {\proclaimfont a canonical basis $\Ga_1,\Ga_2,\Ga_3$
is determined by a marking at~$z$ uniquely
up to conjugation by $\Ga_1\Ga_2\Ga_3$, \ie, up to
the central element $\tau^3\in\BG3$.}

A canonical basis defines an isomorphism $\rho_z\colon G_z\to G$ to the
`standard' free group $G=\<\Ga_1,\Ga_2,\Ga_3>$. This isomorphism
is determined by a marking at~$z$ up to~$\tau^3$.
Below,
{\it all braids involving~$\rho_z$ are considered
up to a power of~$\tau^3$}.
The isomorphism~$\rho_z'$ defined by
the cyclic permutation $e_2,e_3,e_1$ of the bold edges is
given by $\rho_z'=\tau\circ\rho_z$.

\remark
In a similar way, one can define a canonical basis and isomorphism
$\rho_z\colon G_z\to G$ for a nonsingular \white-vertex~$z$. The
basis and the isomorphism are determined up to a power of~$\tau^3$
by an ordering of the two \bold- edges incident to~$z$. Our choice
of \black-vertices is motivated by the fact that we will apply the
results to skeletons.
\endremark

\subsection{Assumptions and settings}\label{assumptions}
For the rest of this section,
we make the following assumptions about~$\Gamma$:
\roster
\item\local{generic}
$\Gamma$ has no monochrome vertices, all its \black-vertices have
valency $\le6$, and all its \white-vertices have valency $\le4$;
\item\local{connected}
the union of all \black-- and \white-vertices of~$\Gamma$ and its
\bold- edges is connected;
\item\local{base}
$\Gamma$ has at least one nonsingular \black-vertex.
\endroster
Note that Condition~\loccit{generic}
means that the curve is
generic within its fiberwise deformation class, and
\loccit{connected} can be satisfied after a sequence of
equisingular perturbations and degenerations,
\cf.~\cite{DIK.elliptic}.
Thus, the only true restriction is~\loccit{base}.
In particular,
any maximal dessin satisfies~\loccit{generic}
and~\loccit{connected}, and the remaining
Condition~\loccit{base} rules out four series of
maximal curves: those whose skeleton is a simple cycle (one curve
in~$\Sigma_k$ for each $k\ge1$)
or
a linear tree (two curves in~$\Sigma_1$ and three curves
in~$\Sigma_k$ for $k\ge2$; a curve is determined by the number of
hanging edges in the skeleton). All these curves are irreducible.

Chose and fix
the `fiber at infinity'~$F_\infty$ over a
point $\infty\notin\Gamma$ and a proper section~$S$,
see~\ref{proper.section}. Denote by $S^\circ\subset\Cp1\cong S^2$
the affine plane $\Cp1\sminus\infty$ punctured at the singular
fibers of~$B$. (Since $S$ is a section,
$S^\circ$ can as well be regarded as a subset of~$S$.)

As above, let $G=\<\Ga_1,\Ga_2,\Ga_3>$ be the free group on three
generators.
Fix a marking of~$\Gamma$ and consider the corresponding
isomorphisms $\rho_z\colon G_z\to G$, see~\ref{basis}.
Given a path~$\Gg$ in~$S^\circ$ connecting two nonsingular
\black-vertices~$z'$ and~$z''$, consider the monodromy
$\tilde m_\Gg\colon G_{z'}\to G_{z''}$ and define the automorphism
$m_\Gg=\rho_{z''}\circ\tilde m_\Gg\circ\rho_{z'}^{-1}$ of~$G$.
It is a braid (due to the fact that $S$ is proper).
We consider $m_\Gg$ as an element of the reduced
group $\BG3/\tau^3$, thus removing the ambiguity in the definition
of~$\rho$. In the special case $z'=z''$, \ie, when $\Gg$ is a
loop, $m_\Gg$ is a well defined element of~$\BG3$. It can be
recovered from its image in $\BG3/\tau^3$ using the following
obvious statement.

\proposition\label{prop.degree}
The degree of the monodromy $\tilde m_\Gg\colon G_z\to G_z$ defined by
a simple loop $\Gg$ in~$S^\circ$
is equal to
the total multiplicity $\sum\mult F_i$
\rom(see Table~\ref{table.j}\rom)
of the singular fibers of~$B$ encompassed by~$\Gg$, \ie, separated
by~$\Gg$ from~$\infty$.
\qed
\endproposition

\subsection{The monodromy}\label{monodromy}
To uniformize the formulas below, we use the convention
$e_{3t+i}=e_i$, $i=1,2,3$, $l\in\Z$,
for the ordered edges incident to a given
nonsingular \black-vertex
(\cf. similar convention for the braid group in~\ref{B3}).

Let $z'$, $z''$ be two nonsingular \black-vertices, connected
by the path~$\Gg$ in~$\Gamma$ formed by two \bold- edges incident to
the same \white-vertex. Denote $m_\Gg=\bb_{i,j}\in\BG3/\tau^3$,
where $i$,~$j$ are the indices of the edges
constituting~$\Gg$ at~$z'$ and~$z''$, respectively. Then
\begin{equation}
\bb_{i,i+1}=\Gs_i,\quad \bb_{i+1,i}=\Gs_{i}^{-1},
 \quad\text{and}\quad
 \bb_{i,i}=\Gs_i\Gs_{i-1}\Gs_i.\label{eq.bb}
\end{equation}
More generally, let $s\ge0$ be an integer, and
let~$\Gg$ be a simple path from~$z'$ to~$z''$ composed
of $2s$ \bold- edges, $(s+1)$ \white-vertices, and $s$
\black-vertices of valency~$4$. Perturb~$\Gg$ so
that each singular \black-vertex
is circumvented in the counterclockwise
direction, and denote by $\bb_{i,j}\ix{s}\in\BG3/\tau^3$
the resulting monodromy.
Then, for all integers $s,t\ge0$, there is a
reciprocity relation
\begin{equation}
\bb_{j+1,i}\ix{t}\cdot\bb_{i+1,j}\ix{s}\cdot\Gs_{i}^{s+t+2}=1,
\label{eq.recipr}
\end{equation}
which can be used to find $\bb_{*,*}\ix{s}$ in terms of
$\bb_{*,*}\ix{0}=\bb_{*,*}$.
One has
$$
\bb_{i,i+1}\ix{s}=\Gs_{i+1}^{-s}\Gs_i,\quad
 \bb_{i+1,i}\ix{s}=\Gs_{i}^{-s-1},
 \quad\text{and}\quad
 \bb_{i,i}\ix{s}=\Gs_i^{-s-2}\Gs_{i+1}^{-1}.
$$

Now, let $\Gg$ be the loop composed of a small counterclockwise
circle around a \cross-vertex of valency~$2d$
connected along a \solid-
edge~$e_i'$ (see Figure~\ref{fig.basis}, left) to a nonsingular
\black-vertex~$z$. The resulting monodromy $\bc_i\ix{d}=m_\Gg\in\BG3$ is
given by
\begin{equation}
\bc_i\ix{d}=\Gs_{i+1}^d.\label{eq.bc}
\end{equation}

Finally, consider a chain of distinct \bold- edges starting from
an edge~$e_i$
at a
nonsingular \black-vertex~$z$ and ending at a singular vertex. Let
$\Gg$ be a simple loop at~$z$ encompassing all vertices of the
chain (except $z$ itself) and oriented in the
counterclockwise direction, and let $\bl_i\ix{d}=m_\Gg\in\BG3$ be
the monodromy, where $d=\deg\bl_i\ix{d}$. (If the chain contains
$s$ \black-vertices of valency~$4$, then $d$ can take the values
$4s$, $4s+2$, or $4s+3$, depending on whether the chain ends at a
\black-vertex of valency~$4$, \white-vertex, or \black-vertex of
valency~$2$.) One has
\begin{equation}
\bl_i\ix{4s}=\Gs_i^{-s}\Gs_{i-1}^{-s}\tau^{3s}
 \quad\text{and}\quad
 \bl_i\ix{4s+\epsilon}=\Gs_i^{-2s-5+\epsilon}\Gs_{i+1}^{-1}\tau^{3s+3},
\label{eq.bl}
\end{equation}
where $\epsilon=2$ or~$3$.

\subsection{Proofs}
But for the choice of the trivialization of the ruling, which is
also accountable for the $\tau^3$-ambiguity, the monodromy~$m_\Gg$
is local with respect to~$\Gg$, and it can be found using the
description of the geometry of the fibers given
in~\ref{tg->topology}. We do use this straightforward approach to
establish relations~\eqref{eq.bb} and~\eqref{eq.bc}. The
expression for~$\bl_i\ix{4s}$ in~\eqref{eq.bl} follows from
Proposition~\ref{prop.degree} and the obvious relation
$$
\bl_i\ix{4s}=\bb_{j,i}\ix{s}\cdot\bb_{i,j}\ix{s},\quad j\in\Z,
$$
in $\BG3/\tau^3$,
which is due to our convention that the paths are perturbed so as
to circumvent all singular vertices in the counterclockwise
direction.

For the rest, we observe that the monodromy related to a fragment
of~$\Gamma$ can be found in any other dessin containing this
fragment. The reciprocity relation~\eqref{eq.recipr} is obtained
assuming that the two paths resulting in the two $\bb_{*,*}$
monodromies form the boundary (oriented in the clockwise
direction) of a single region~$R$ of the skeleton of the dessin,
so that $R$ contains a single \cross-vertex. (The factor
$\Gs_i^{s+t+2}$
in the relation is, in fact, $\bc_{i-1}(s+t+2)$.) The expressions
for~$\bl_i$ are obtained in a similar way: we close the unused
\bold- edges $e_{i-1}$, $e_{i+1}$ at~$z$ `around' the chain of
edges in question and place a single \cross-vertex at the center
of the resulting region~$R$. Computing the monodromy
around~$\partial R$ gives the relations
$$
\bc_{i-1}\ix{2s+5-\epsilon}\cdot\bl_i\ix{4s+\epsilon}=
 \bl_i\ix{4s+\epsilon}\cdot\bc_{i+1}\ix{2s+5-\epsilon}=
 \bb_{i-1,i+1}
$$
in $\BG3/\tau^3$
(where $\epsilon=2$ or~$3$), which can be used to find~$\bl_i$.

\corollary\label{reducible.curve}
A maximal trigonal curve~$B$ is reducible if and only if
all vertices of its skeleton~$\Sk$
are nonsingular \rom(\ie, have valency~$3$\rom) and $\Sk$
admits a marking with the following properties\rom:
\roster
\item\local{1}
each hanging edge has index~$1$ at the \rom(only\rom)
vertex incident to it\rom;
\item\local{2}
any other edge has indices $(1,1)$, $(2,3)$, or $(3,2)$ at
its two endpoints.
\endroster
\endcorollary

\remark
Clearly, a marking at any vertex of~$\Sk$ extends to at most one
marking satisfying Condition~\iref{reducible.curve}2.
If $\Sk$ has a
hanging edge, it admits at most one marking
satisfying~\iref{reducible.curve}1 and~\ditto2.
\endremark

\remark
Corollary~\ref{reducible.curve} still makes sense for a trigonal
curve~$B$, not necessarily maximal, whose dessin~$\Gamma$
satisfies Conditions~\iref{assumptions}{generic}
and~\ditto{connected}. In this case, the existence of a marking as
in Corollary~\ref{reducible.curve} is necessary for~$B$ to be
reducible; in general, it is not sufficient.
\endremark

\proof
Let $B^\circ$ be the portion of the curve over~$S^\circ$, and let
$\pr\colon B^\circ\to S^\circ$ be the restriction of the projection
$\Sigma_k\to\Cp1$. It is a triple
covering whose monodromy is obtained by downgrading
the braid monodromy to the symmetric
group~$\SG3$. From~\eqref{eq.bl} it follows that the monodromy
about a singular \black-vertex acts transitively on the decks
of~$\pr$, and hence any curve with such a vertex is irreducible.
(As this argument is local, it applies as well to the four
exceptional series mentioned in~\ref{assumptions}, proving that
they are all irreducible.)

Assume that $B$ is reducible. Then it contains as a component
a section of the ruling. Any such section $B_1\subset B$ defines a
marking of~$\Sk$:
one assigns index~$1$ to the point $B_1\cap F_z\in B_z$,
see~\ref{basis}, and
Conditions~\iref{reducible.curve}1 and~\ditto2 merely list all
monodromies~$\bl_i(3)$ and~$\bb_{i,j}$ preserving~$p_1$.
Conversely, for any marking as in the statement, the points
$p_1\in B_z$ over all \black-vertices $z\in\Gamma(B)$ belong to a
deck of~$\pr$ which is preserved by the monodromy. (One needs to
take into account the obvious fact that, for a maximal curve
without singular \black-vertices, the
inclusion homomorphism $\pi_1(\Sk)\to\pi_1(S^\circ)$ is an
isomorphism.) Hence, the curve contains a section of the ruling as
a component.
\endproof

\section{The construction}\label{S.construction}

Proofs of Theorems~\ref{th.main} and~\ref{th.reducible} are based
on the existence of large Zariski $k$-plets of maximal trigonal
curves in Hirzebruch surfaces.

\proposition\label{C.curves}
For each integer $k\ge2$, there exists a collection of
$$
C(k-1)=\frac1k\binom{2k-2}{k-1}
$$
pairwise distinct fiberwise deformation
families of irreducible maximal
trigonal curves $B\subset\Sigma_k$ with the following
properties\rom:
\roster
\item
each curve has one fiber of type~$\tilde\bA_0^{**}$, one
fiber of type $\tilde\bA_{5k-3}$, and $k$ fibers
of type~$\tilde\bA_0^*$ \rom(and no other singular fibers\rom)\rom;
\item
none of the curves admits a fiberwise self-deformation inducing a
non-trivial permutation of the singular fibers of the curve.
\endroster
\endproposition

\begin{figure}[ht]
\centerline{\epsffile{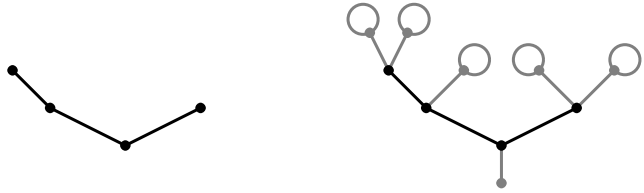}}
\caption{Extending a binary tree~$T$ to a skeleton~$\Sk(T)$}\label{fig.tree}
\end{figure}

\proof
Denote by~$\CT_s$, $s\ge1$, the set of all binary rooted trees on
$s$~vertices. Recall that the cardinality of~$\CT_s$ is
given by the Catalan number~$C(s)$,
$$
\#\CT_s=C(s)=\frac1{s+1}\binom{2s}s.
$$
Each tree $T\in\CT_s$ admits a standard `monotonous' geometric
realization $\ls|T|\subset\R^2$, see Figure~\ref{fig.tree}, left.
For example, one can map the level~$l$, $l\ge0$, vertices
of~$T$ to the points
$v_{l,i}=(-1+(2i+1)/2^l,l)$, $i=0,\ldots,2^l-1$, so that the left
(respectively, right) edge originating at~$v_{l,i}$
connects~$v_{l,i}$ to $v_{l+1,2i}$
(respectively, $v_{l+1,2i+1}$).

Pick a tree $T\in\CT_{k-1}$ and extend its
geometric realization $\ls|T|\subset\R^2\subset\Cp1$
to a skeleton~$\Sk(T)$
as follows: mark the root of~$T$
by
adding a monovalent vertex at $(0,-1)$ and connecting it
to~$v_{0,0}$ by an edge,
and complete the valency of each vertex of~$\ls|T|$ to
three by
replacing the missing branches with
`leaves', each leaf
consisting of a vertex (at an appropriate point $v_{l,i}$, $l>0$),
a loop at this vertex, and a stem connecting the vertex to the
point $v_{l-1,[i/2]}$.
(See Figure~\ref{fig.tree}, right, where the trunk and
the $k$~leaves added to~$\ls|T|$ are shown
in grey.)

The resulting skeleton~$\Sk(T)$ has one monovalent and $(2k-1)$
trivalent vertices; its faces are $k$~monogons (the interiors of the
leaves) and one $(5k-2)$-gon (the outer region). Furthermore, one
can easily observe that
none of $\Sk(T)$ has a nontrivial automorphism and that
two skeletons $\Sk(T_1)$, $\Sk(T_2)$ are isomorphic
if and only if
$T_1=T_2$ in~$\CT_{k-1}$.
Here, the key observation is the fact that the root of the
original tree~$T$
is `marked' by the only monovalent vertex of the skeleton~$\Sk(T)$.
Hence,
any isomorphism of the skeletons would induce an isomorphism of
oriented \emph{rooted} trees (as it also preserves the orientation
of~$\base$). In particular, essentially by
its very definition, an oriented \emph{rooted} tree never
admits an orientation
preserving automorphism.

Applying Theorem~\ref{B->Sk},
one obtains $\#\CT_{k-1}=C(k-1)$ deformation families of
maximal trigonal curves with the desired properties.
(Each curve is irreducible since it has a type~$\tilde\bA_0^{**}$
singular fiber.)
\endproof

\proposition\label{C.curves.real}
If $k=2s$ is even, then $C(s-1)$ of the trigonal curves given by
Proposition~\ref{C.curves} are real \rom(with respect to some real
structure on~$\Sigma_k$\rom). All other curves \rom(and all curves
for $k$ odd\rom) split into pairs of complex conjugate curves.
\endproposition

\proof
A maximal trigonal curve is real if and only if its skeleton is
symmetric with respect to some \emph{orientation reversing}
involution of the base~$\base$
(\cf.~\cite{DIK.elliptic},
\S5 and
especially Corollary~5.1.8,
where real dessins of real
curves are considered: clearly, a symmetric skeleton can be
completed to a symmetric dessin, due to the results
of~\cite{DIK.elliptic} cited above, a symmetric dessin
gives rise to a real $j$-invariant, and the further passage from
the $j$-invariant to a trigonal curve is equivariant, \cf. the
proof of Proposition~\ref{j->B}). A binary rooted tree can be
symmetric only if its number of vertices is odd, and all symmetric
trees in $\CT_{2s-1}$ can be parametrized by their `left halves',
\ie, by $\CT_{s-1}$.
\endproof

\subsection{The braid monodromy}\label{C.monodromy}
In this section, we apply the results of~\ref{monodromy} to
describe the braid monodromy of the curves given by
Proposition~\ref{C.curves}.

Fix a curve~$B$ corresponding to a tree $T\in\CT_{k-1}$
and let $\Gamma=\Gamma(B)$, $\Sk=\Sk(B)$.
Let $v_{0,0}$ be the root of the original tree~$T$.
Denote by $\Gamma_{\!\CROSS}$ the set of \cross-vertices
of~$\Gamma$ of valency~$2$ (equivalently, the set of
type~$\tilde\bA_0^*$ singular fibers of~$B$).
Each vertex $u\in\Gamma_{\!\CROSS}$ can be encoded by a
word~$w_u$ in the alphabet $\{r,l\}$ as follows: let $\bar u$ be the
\black-vertex in the leaf encompassing~$u$, and
let $\xi_u$ be the simple path in~$\Sk$ from~$v_{0,0}$
to~$\bar u$; starting
from~$v_{0,0}$ and the empty word,
walk along $\xi_u$ and, at each
vertex, add to the word~$r$ or~$l$ if the right (respectively,
left) branch is chosen at this vertex. (For example, in
Figure~\ref{fig.tree}, the \cross-vertices encompassed by the five
leaves are encoded, from right to left, by the words $rr$, $rl$,
$lr$, $llr$, and $lll$.) Order $\Gamma_{\!\CROSS}$
lexicographically, with $r<l$. (This is the right
to left order in the standard geometric realization of the graph,
\cf. Figure~\ref{fig.tree}.)

As in~\ref{assumptions}, pick a point $\infty\in\Cp1\sminus\Gamma$
and denote by~$S^\circ$ the plane $\Cp1\sminus\infty$ punctured at
the singular vertices of~$\Gamma$. Take~$v_{0,0}$ for the base
point, and consider the basis $\Gg_u$, $u\in\Gamma_{\!\CROSS}$,
$\Gd_{\!\CROSS}$, $\Gd_{\BLACK}$ for $\pi_1(S^\circ,v_{0,0})$
defined as follows:
\roster
\item
$\Gg_u$, $u\in\Gamma_{\!\CROSS}$, is the loop in~$\Sk$ formed by
the circumference of the leaf surrounding~$u$ connected
to~$v_{0,0}$ by the simple path~$\xi_u$;
\item
$\Gd_{\!\CROSS}$ is a small circle surrounding the \cross-vertex of
valency $10k-4$, connected to~$v_{0,0}$ by the left \solid- edge
at~$v_{0,0}$;
\item
$\Gd_{\BLACK}$ is a small circle surrounding the singular
\black-vertex, connected to~$v_{0,0}$ by the \bold- edge.
\endroster
(All loops are oriented in the counterclockwise direction.) Then,
the braid monodromy $\pi_1(S^\circ,v_{0,0})\to\BG3$ is given by
the following relations:
$$
\Gg_u\mapsto\bar w_u\Gs_3\bar w_u^{-1},\quad
 \Gd_{\!\CROSS}\mapsto\Gs_1^{5k-2},\quad
 \Gd_{\BLACK}\mapsto\Gs_1\Gs_2,
$$
where $\bar w_u$ is the braid obtained from~$w_u$ by replacing
each instance of~$r$ and~$l$ with $\Gs_2$ and $\Gs_1^{-1}$,
respectively.

\proof
We choose a marking at each nonsingular \black-vertex of~$\Sk$
so that $e_2$ is the edge pointing downwards (and hence $e_1$
and~$e_3$ are, respectively, the left and right branches of the
tree). Then $\Gd_{\!\CROSS}\mapsto\bc_3(5k-2)$, see~\eqref{eq.bc},
$\Gd_{\BLACK}\mapsto\bl_2(2)$, see~\eqref{eq.bl}, and the image
of each element~$\Gg_u$ is found by composing appropriate
monodromies $\bb_{i,j}$, see~\eqref{eq.bb}.
\endproof

\subsection{Proof of Theorem~\ref{th.main}}\label{proof.main}
Let~$k$ and~$\epsilon$ be as in the statement. Note that $k\ge3$
and $\epsilon\le[k/2]$.
Pick one of the
trigonal curves $B\subset\Sigma_k$ given by
Proposition~\ref{C.curves}. In order to convert~$B$ to a plane
curve, we need to perform $(k-1)$ elementary
transformations.
We choose the transformations so as to contract the
type~$\tilde\bA_0^{**}$ fiber of~$B$, $[k/2]$ of its $k$
type~$\tilde\bA_0^*$ fibers, and $[(k-3)/2]$ nonsingular fibers.
In the type~$\tilde\bA_0^{**}$ fiber and $\epsilon$
type~$\tilde\bA_0^*$ fibers
the blow-up centers are chosen outside of the
curve and the exceptional section; in each other fiber the blow-up
center is taken on a branch of~$B$ transversal to
the fiber.
The total number of deformation families
thus obtained is
$$
Z(m)=C(k-1)\cdot\binom{k}{[k/2]}\cdot\binom{[k/2]}\epsilon,
$$
the three factors standing, respectively, for the choice of~$B$,
the choice of $[k/2]$ of its $k$ type~$\tilde\bA_0^*$ fibers to be
contracted, and the choice of~$\epsilon$ of the $[k/2]$ fibers where
the blow-up center is not on~$B$. In each case, the transform is an
irreducible curve
$\tilde C=B_1\subset\Sigma_1$ with the
set of singularities
\begin{equation}
\bA_{5k-3}+\bE_6+\epsilon\bD_5+
 \Bigl(\Bigl[\frac{k}2\Bigr]-\epsilon\Bigr)\bA_2+
 \Bigl[\frac{k-3}2\Bigr]\bA_1,\label{eq.proof}
\end{equation}
so that all points except the first $\bA_{5k-3}$ are in the
exceptional section $E_1\subset\Sigma_1$ and the local intersection index
of~$\tilde C$ and~$E_1$ at each singular point is minimal possible
(\ie, $2$ at a double point and~$3$ at a triple point).
Blowing~$E_1$ down, one obtains an irreducible plane curve~$C$ of
degree $2k+2+\epsilon=m$. Since the combinatorial data of~$C$ are
determined by the those of $\tilde C+E_1$, all curves
thus obtained share the same set of singularities.

The fundamental groups $\pi_1(\Cp2\sminus C)$ are all abelian
due to Proposition~\ref{abelian.group}.
\qed

\subsection{Proof of Theorem~\ref{th.real}}\label{proof.real}
A real curve~$C$ is obtained from a real trigonal curve
$B\subset\Sigma_k$; hence, $k=2s$ is even and the number of
real trigonal curves is given by Proposition~\ref{C.curves.real}.
Next, one should choose a real (\ie, invariant under the complex
conjugation) collection of blow-up centers for the elementary
transformations converting~$B$ to~$\tilde C$,
see~\ref{proof.main}. Since the $k$ type~$\tilde\bA_0^*$ singular
fibers of~$B$ split into $k/2=s$ conjugate pairs and $k/2=s$
blow-up centers should be chosen in these fibers, $s$ must also be
even, $s=2t$, and the number of choices is $\binom st$: one
chooses $t$ of the $s$~conjugate pairs. Finally, $\epsilon=0$ as
one cannot choose only one special fiber with the blow-up center
not on the curve: for the transformation to be real, the conjugate
fiber would have to have the same property.
\qed

\remark
It is worth mentioning that, in the settings of
Theorem~\ref{th.real}, each deformation class containing a real
curve splits into at least $t^2$ equisingular real deformation
classes. Indeed, let $\Rp1\subset\Cp1$ be the real part of the
base of the ruling. It contains the singular \black-vertex
of~$\Gamma$, the root of the original tree, and the \cross-vertex
of~$\Gamma$ of valency $10k-4$. Thus, the singular fibers of~$B$
divide~$\Rp1$ into two distinguishable intervals, and each of
the $2t-2$
nonsingular fibers containing blow-up centers can be chosen either
in a conjugate pair or over one of the two intervals. The number
of choices is the number of ordered pairs $(a,b)\in\Z\times\Z$
such that $a,b\ge0$, $a+b\le2t-2$, and $a+b$ is even. It is $t^2$.
\endremark

\subsection{Proof of Theorem~\ref{th.reducible}}\label{proof.reducible}
The proof is similar to
that of Theorem~\ref{th.main}.
Let $k=m-5$, and pick
one of the trigonal curves $B\subset\Sigma_k$ given by
Proposition~\ref{C.curves}. Blow up the only singular point of~$B$
and blow down the corresponding fiber. Repeat this procedure
$(k-2)$ times. The result is an irreducible curve
$B_2\subset\Sigma_2$ which intersects the exceptional section at a
nonsingular point~$P$ with multiplicity $(k-2)$ and has a
type~$\bA_{3k+1}$ singular point~$Q$ in the fiber~$F_P$ through~$P$.
Now, add the fiber~$F_P$ as a component, perform an elementary
transformation to contract the type~$\tilde\bA_0^{**}$ fiber
of~$B_2$ to a type~$\bE_6$ singular point in the exceptional
section, and blow down the exceptional section. (The fiber~$F_P$
is added to the curve to make sure that, during the deformations,
the intersection point~$P$ and the singular point~$Q$ remain in
the same fiber.)
The result is a
plane curve~$C$ of degree~$m$. Clearly, all $C(k-1)$ curves
obtained in this way share the same combinatorial data. The
fundamental groups $\pi_1(\Cp2\sminus C)$ are all abelian due to
Proposition~\ref{abelian.group.ex}.

The count for the number of real curves is based on
Proposition~\ref{C.curves.real}: from the construction it follows
that a family contains a real curve if and only if the original
trigonal curve~$B$ is real.
\qed

\subsection{A remark on deformations}\label{rem.deformations}
From the construction (creating a branch of type~$\bE_{12}$ and,
in~\ref{proof.reducible}, adding a linear component)
it follows that any equisingular deformation
of the plane curve~$C$ must preserve the type~$\tilde\bA_0^{**}$
and type~$\tilde\bA_{5k-3}$ singular fibers of the original
trigonal curve~$B$. Since all other singular fibers of~$B$ are of
type~$\tilde\bA_0^*$ and $B$ is maximal, the resulting deformation
of~$B$ is fiberwise, see Proposition~\ref{no.degeneration}. The
blow-up centers chosen in the branches of~$B$ transversal to its
type~$\tilde\bA_0^*$ singular fibers (see~\ref{proof.main}) should
stay fixed, as otherwise the type of singularity of~$C$ at~$O$
would change. (This observation is also crucial in the proof of
Theorem~\ref{alt.main} below.) In~\ref{proof.main},
a blow-up center in a nonsingular
fiber of~$B$ may move to a type~$\tilde\bA_0^*$ singular fiber
(not containing another blow-up center), to the branch of~$B$
\emph{tangent}
to the fiber. This degeneration corresponds to one of the
branches of one of the type~$\bA_1$ points of~$\tilde C$,
see~\eqref{eq.proof}, becoming tangent to the fiber; it is
equisingular for~$C$. Clearly, these modifications do not affect
the number of deformation families.

\section{Further applications}\label{S.applications}

In this section, we present a slight modification of the
construction used in Proposition~\ref{C.curves} and discuss
a few further applications.

\proposition\label{alt.curves}
For each integer $k\ge2$, there exists a collection of
$C(k-1)$
pairwise distinct fiberwise deformation
families of pairs $(B,F)$, where $B\subset\Sigma_k$ is an
irreducible maximal
trigonal curve with
one fiber of type~$\tilde\bA_{5k-2}$ and $(k+1)$ fibers
of type~$\tilde\bA_0^*$ \rom(and no other singular fibers\rom) and
$F$ is a distinguished type~$\tilde\bA_0^*$ fiber of~$B$.
None of the curves admits a fiberwise self-deformation inducing a
non-trivial permutation of its singular fibers preserving~$F$.
\endproposition

\proof
Modify the construction of
Proposition~\ref{C.curves} by replacing the monovalent vertex with
an extra leaf attached to the root of
the original tree~$T$ and selecting the
corresponding type~$\tilde\bA_0^*$ fiber for~$F$. All curves
obtained are irreducible due to Corollary~\ref{reducible.curve}:
to show that a marking as in the corollary does not exist, it
suffices to consider the two leaves attached to any maximal (in
the partial order defined by level) vertex of~$T$.
\endproof

\theorem[Rigid plane curves]\label{alt.main}
For each odd integer $m=2k+1\ge5$,
there is a set of singularities shared by
$$
\Za(m)\ge\frac12\binom{2k-2}{k-1}
$$
pairwise distinct equisingular deformation families of
irreducible plane curves~$C_i$ of
degree~$m$. Within each family, all curves are projectively
equivalent and defined over an algebraic number field.
\endtheorem

\remark
The set of singularities constructed in the proof has a point of
type~$\bA_{5k-2}$ and a point of transversal intersection of
$(k-1)$ branches of type~$\bA_4$. One has $\Za(5)=1$, and the only
curve of degree~$5$ given by the theorem is the well known quintic
with the set of singularities $\bA_8+\bA_4$, see~\cite{quintics};
it is defined over~$\Q$.
(Note that in this case the fundamental group
$\pi_1(\Cp2\sminus C)$ is abelian, see~\cite{groups}.)
For large values of~$m$, the count
$\Za(m)$ grows faster than $a^m$ for any $a<2$.
\endremark

\remark
The curves given by Theorem~\ref{alt.main} seem to be good
candidates for examples of exponentially large \emph{arithmetic
Zariski $k$-plets} in the sense of Shimada, see~\cite{Shimada1},
\cite{Shimada2}. At present, I do not know whether all/some of the
curves~$C_i$ are indeed Galois conjugate over an algebraic number
field (except the trivial case of pairs of complex conjugate
curves). Whether the pairs $(\Cp2,C_i)$ or complements
$\Cp2\sminus C_i$ are homeomorphic is also an open question.
\endremark

\proof
Similar to~\ref{proof.main}, we start with a trigonal
curve~$B\subset\Sigma_k$ as in
Proposition~\ref{alt.curves}, perform $(k-1)$ elementary
transformations to convert~$\Sigma_k$ to~$\Sigma_1$, and blow down
the exceptional section of~$\Sigma_1$ to get a plane curve. The
$(k-1)$ blow-up centers are taken in type~$\tilde\bA_0^*$
singular fibers of~$B$, on the branch of~$B$ transversal to the
fiber. (This choice makes the construction rigid, so that the
resulting plane curves have $0$-dimensional moduli spaces and are
defined over algebraic number fields.
Indeed, since $B$ itself is defined over a certain algebraic
number field~$\Bbbk$, see remark after~\ref{def.maximal}, all
its singular fibers~$F_j$ are defined over a finite extension
of~$\Bbbk$, and so are the intersection points $B\cap F_j$. Hence,
each curve~$C_i$
is also defined over a finite extension of~$\Bbbk$.)
The total number of choices
is $C(k-1)$ (for the pair $(B,F)$\,) times $k(k+1)/2$ (for the
choice of $(k-1)$ singular fibers containing the blow-up centers).
Since, in each
skeleton, the distinguished leaf can be chosen in $(k+1)$ ways, we
divide the resulting count by $(k+1)$.
\endproof

\subsection{Elliptic surfaces}\label{surfaces}
Below, an \emph{elliptic surface} is a compact complex
surface~$X$ with a
distinguished rational pencil of elliptic curves,
\ie, elliptic fibration over a rational base.
We assume that the pencil has no multiple fibers; then it is
unique unless the topological Euler characteristic of~$X$
is~$24$, \ie, $X$ is a $K3$-surface.
By a
\emph{fiberwise deformation} of elliptic surfaces we mean a
deformation preserving the elliptic pencil and the types of its
singular fibers. All surfaces mentioned in
Theorems~\ref{th.surfaces} and~\ref{alt.surfaces} are defined over
algebraic number fields.

\theorem\label{th.surfaces}
For each integer $s\ge1$, there are $C(2s-1)$ distinct fiberwise
deformation families of
Jacobian relatively minimal
elliptic surfaces of topological
Euler characteristic $\chi=12s$ and having
one fiber of type~$\tilde\bA_0^{**}$, one
fiber of type $\tilde\bA_{10s-3}$, and $2s$ fibers
of type~$\tilde\bA_0^*$ \rom(and no other singular fibers\rom).
\endtheorem

\theorem\label{alt.surfaces}
For each integer $s\ge1$, there are
at least
\hbox{$C(2s-1)/(2s+1)$}
distinct fiberwise
deformation families of
Jacobian relatively minimal
elliptic surfaces of topological
Euler characteristic $\chi=12s$ and having
one fiber of type $\tilde\bA_{10s-2}$ and $(2s+1)$ fibers
of type~$\tilde\bA_0^*$ \rom(and no other singular fibers\rom).
\endtheorem

\proof[Proof of Theorems~\ref{th.surfaces} and~\ref{alt.surfaces}]
The statements follow from Propositions~\ref{C.curves}
and~\ref{alt.curves} applied to $k=2s$. Each surface is obtained
as the minimal resolution of singularities of the double covering
of~$\Sigma_k$ branched over the exceptional section~$E$ and a
trigonal curve~$B$ given by the appropriate proposition.
\endproof

\remark
Let~$X$ be one of the surfaces given by Theorem~\ref{th.surfaces}
or~\ref{alt.surfaces}, and let $L=H_2(X)$ be its intersection
lattice. Consider the sublattice $S\subset L$ spanned by the
components of the pull-back of $B\cup E$. Over~$\Q$, it is spanned
by the section of~$X$, its generic fiber, and the exceptional
divisors over the only singular point of~$B$. Hence, $S$ is
nondegenerate. The advantage of Theorem~\ref{alt.surfaces} is the
fact that, in this case, the orthogonal complement $S^\perp$ is an
even \emph{positive definite} lattice of rank~$2s-2$. Given that
positive definite lattices tend to have many isomorphism classes
within the same genus, one can hope to use Shimada's
invariant~\cite{Shimada1} to distinguish the surfaces
topologically.
\endremark


\let\ref\bibitem
\let\endref\relax
\let\publ\relax
\let\publaddr\relax
\def\by#1{\rm #1,}
\def\paper#1{{\it #1\/},}
\def\book#1{{\it #1\/},}
\def\jour#1{\rm #1,}
\def\yr#1{(#1),}
\def\issue#1{no.~#1,}
\def\vol#1{{\bf#1}}
\def\pages#1{#1.}
\def\finalinfo#1{{#1}}


\end{document}